\documentclass[12pt]{amsart}

\usepackage{url}
\usepackage{graphicx}
\usepackage{amsmath}
\usepackage{amscd}
\usepackage{amsfonts}
\usepackage{amssymb}
\usepackage{color}
\usepackage{fullpage}
\usepackage{mathtools}
\usepackage[pagebackref,hypertexnames=false, colorlinks, citecolor=red,linkcolor=blue, urlcolor=red]{hyperref}
\usepackage{todonotes}
\usepackage{enumitem}
\usepackage{comment}

\numberwithin{equation}{section}

\newtheorem{theorem}{Theorem}[section]
\newtheorem{lemma}[theorem]{Lemma}
\newtheorem{proposition}[theorem]{Proposition}

\newtheorem{corollary}[theorem]{Corollary}

\newtheorem*{theorem*}{Theorem}
\newtheorem*{lemma*}{Lemma}

\theoremstyle{definition}
\newtheorem{example}[theorem]{Example}
\newtheorem{remark}[theorem]{Remark}

\newtheorem*{remark*}{Remark}
\newtheorem*{definition*}{Definition}
\newtheorem*{problem*}{Problem}

\newcommand{\be}{\begin{equation}}
\newcommand{\ee}{\end{equation}}
\newcommand{\bes}{\begin{equation*}}
\newcommand{\ees}{\end{equation*}}

\newcommand{\cA}{\mathcal{A}}
\newcommand{\cB}{\mathcal{B}}

\newcommand{\cE}{\mathcal{E}}
\newcommand{\cF}{\mathcal{F}}

\newcommand{\cH}{\mathcal{H}}
\newcommand{\cI}{\mathcal{I}}

\newcommand{\cM}{\mathcal{M}}

\newcommand{\bB}{\mathbb{B}}
\newcommand{\bC}{\mathbb{C}}
\newcommand{\C}{\mathbb{C}}
\newcommand{\bD}{\mathbb{D}}

\newcommand{\bM}{\mathbb{M}}
\newcommand{\bN}{\mathbb{N}}

\newcommand{\ol}{\overline}

\newcommand{\diag}{\operatorname{diag}}

\newcommand{\id}{\operatorname{id}}

\newcommand{\mlt}{\operatorname{Mult}}

\newcommand{\spn}{\operatorname{span}}

\newcommand{\fB}{{\mathfrak{B}}}
\newcommand{\fC}{{\mathfrak{C}}}
\newcommand{\fD}{{\mathfrak{D}}}

\newcommand{\fV}{{\mathfrak{V}}}
\newcommand{\fW}{{\mathfrak{W}}}

\newcommand{\GL}{\mathrm{GL}}

\begin{document}

\title{Isomorphisms and representations of noncommutative function algebras}

\author{Marcel Scherer}
\address{M.S., Department of Mathematics, Technion Israel Institute of Technology\\
 Technion City, Haifa\; 3200003\\
 Israel}

 \author{Orr Moshe Shalit}
 \address{O.M.S., Faculty of Mathematics and The Helen Diller Quantum Center\\
 Technion Israel Institute of Technology\\
 Technion City, Haifa\; 3200003\\
 Israel}
 
 \urladdr{\href{https://oshalit.net.technion.ac.il/}{\url{https://oshalit.net.technion.ac.il/}}}

 \author{Eli Shamovich}
 \address{E.S., Department of Mathematics\\
 Ben-Gurion University of the Negev\\
 Beer-Seva\; 8410501\\
 Israel}

\author{Mansi Suryawanshi}
\address{M.S., Department of Mathematics, Technion Israel Institute of Technology\\
 Technion City, Haifa\; 3200003\\
 Israel}

 \thanks{The work of E. Shamovich is partially supported by BSF grant no.\ 2022235.
 M. Scherer was supported by a Technion Fellowship.
 }
 \subjclass[2020]{47A13, 46L52, 47L55, 46L07}
 \keywords{Noncommutative function algebra, noncommutative operator ball, noncommutative varieties, operator space}

\addcontentsline{toc}{section}{Abstract}

\begin{abstract}
Let $\fV$ be a noncommutative (NC) variety in a NC operator ball. 
We study the dual operator algebra $H^\infty(\fV)$ of bounded NC functions on $\fV$, and the operator algebra $A(\fV)$ obtained as the closure of NC polynomials in $H^\infty(\fV)$. 
We identify the space of weak-$*$ continuous finite-dimensional representations of $H^\infty(\fV)$ as the similarity envelope $\widetilde{\fV}$ of $\fV$, and equip $\widetilde{\fV}$ with a metric which turns $\widetilde{\fV}$ into a complete invariant for $H^\infty(\fV)$, in the sense that $H^\infty(\fV)$ is weak-$*$ continuously isomorphic to $H^\infty(\fW)$ if and only if $\widetilde{\mathfrak{V}}$ and $\widetilde{\mathfrak{W}}$ are NC biholomorphic via a bi-Lipschitz map. 
Similar results are obtained in the completely bounded category as well as for the algebras $A(\fV)$ when $\fV$ is assumed homogeneous. 
In the case of homogeneous varieties $\fV$ and $\fW$, we obtain a sharper classification result, by which there is a weak-$*$ homeomorphic and (completely) bounded isomorphism from $H^\infty(\fV)$ to $H^\infty(\fW)$ if and only if there is a linear map which gives a bi-Lipschitz bijection between $\widetilde{\fW}$ and $\widetilde{\fV}$. 
The reduction to the linear case in the generality we work in depends on showing that the existence of maps implies the existence of origin-fixing maps, which we show using a novel ``two-circle trick" of independent interest. 
In the important case where the varieties are the full NC operator balls we find a striking rigidity phenomenon: if $A(\bB_E)$ and $A(\bB_F)$ are completely boundedly isomorphic then $E$ and $F$ are completely isometric. 
We also present a variety of examples illustrating the results and their sharpness, showing how several results from the literature fail to generalize to the setting studied here. 
\end{abstract}

\maketitle


\maketitle

\section{Introduction}


A pervasive theme in operator algebra since the early development of the field has been the construction of operator algebras that encode a particular mathematical object, and the subsequent study of the relationship between the algebra and the underlying structure. Examples include the von Neumann algebra of a group or of a measurable action on a space, full and reduced group C*-algebras, and graph C*-algebras. 
The quintessential example is the algebra $C(X)$ of continuous functions on a compact Hausdorff space. 
By Gelfand theory, we know that $C(X)$ and $C(Y)$ are isomorphic as algebras if and only if $X$ is homeomorphic to $Y$. 
In fact, the isomorphisms $C(X) \to C(Y)$ are precisely the maps $f \mapsto f \circ h$ induced by homeomorphisms $h: Y \to X$. This, in itself, is a beautiful example of a categorical equivalence, but combined with Gelfand's characterization of commutative C*-algebras, it hints at something deeper: {\em every operator algebra can be viewed as an algebra of ``functions" on some ``space"}. 
This opens the door to the paradigms of noncommutative analysis and noncommutative geometry.

Our research fits in this tradition of studying operator algebras through an underlying geometric structure. 
The classes of algebras that we study in this paper consist of bounded noncommutative functions on subvarieties of operator balls. 
The underlying geometric structures are, interestingly, {\em not} the subvarieties themselves, but rather their similarity envelopes, endowed with NC holomorphic structure and a certain metric. 
Concretely, the main goal of this paper is to classify algebras of bounded NC functions on NC subvarieties of operator balls via their spaces of finite-dimensional representations, which are concretely identified as certain NC varieties. 

The study of the isomorphism problem for algebras of NC functions is a worthwhile pursuit for several reasons. 
First, it enriches the theory of operator algebras by establishing an intimate familiarity with a large class of nonselfadjoint operator algebras with tractable invariants allowing to compare and distinguish between pairs. 
Second, it provides a common framework for classification programs for classes of nonselfadjoint operator algebras \cite{DRS1,DRS2,KakariadisShalit,SS09}. 
Third, representing previously studied operator algebras as algebras of NC functions gives a handle on them and leads to insights and new results. 
Fourth, it leads to function-theoretic questions and results of independent interest \cite{BelSham,kumarsampat2025zeros,sampat2025weak,sampat2026cyclicity,S}. 
Fifth, it has motivated new developments in more general operator-algebraic settings \cite{barik2026algebras}.
Recently, we have observed that these algebras provide concrete representations of universal operator algebras generated by an operator space and its quotients \cite{SampatShalit25,sampat2025weak}, and this suggests deep interactions between operator space theory and NC analytic function theory.

To briefly state the main results of this paper, let us set some notation (which will be explained in detail in the next section). 
Let $E$ and $F$ be finite-dimensional operator spaces, and let  $\fV$ and $\fW$ be subvarieties of the operator balls $\bB_E$ and $\bB_F$, respectively. 
Let $H^\infty(\fV)$ be the algebra of bounded NC functions on $\fV$, let $A(\fV)$ be the closure of the polynomials in $H^\infty(\fV)$, and similarly for $\fW$. 
We seek to determine necessary and sufficient conditions for the algebras $H^\infty(\fV)$ and $H^\infty(\fW)$ to be (completely) boundedly isomorphic (and likewise for the algebras $A(\fV)$ and $A(\fW)$). 

In \cite{SSS1} the classification problem up to completely isometric isomorphism was solved in the case where $\bB_E = \bB_F = \fB_d$ were both the row ball. 
In \cite{SSS2} the problem was solved for $\bB_E = \bB_F = \fB_d$ when classification was done up to bounded and completely bounded isomorphism. 
The papers \cite{SampatShalit25,sampat2025weak} treated the case where the varieties $\fV,\fW$ are subvarieties of arbitrary, perhaps different operator space balls $\bB_E$ and $\bB_F$, and classification was carried out up to completely isometric isomorphism. 
The best results were obtained when the operator spaces $E$ and $F$ were assumed injective and the varieties were matrix spanning homogeneous varieties. 

In this paper, we consider the isomorphism problem up to bounded and completely bounded isomorphism, where the varieties are subvarieties of arbitrary operator space balls. 
We shall now describe the results for $A(\fV)$ and $H^\infty(\fV)$ in tandem. 
We will refer to the algebras of the type $A(\fV)$ as algebras of uniformly continuous NC functions, and to those of type $H^\infty(\fV)$ as algebras of bounded NC functions, keeping in mind that in the uniformly continuous case we need an additional assumption that the varieties are homogeneous, and that in the bounded case we shall take {\em isomorphism} to mean a weak-$*$ homeomorphic isomorphism. 

In Section \ref{sec:fin_reps}, we obtain an identification of the space of bounded finite-dimensional representations of $A(\fV)$ with the similarity $\widetilde{\ol{\fV}}$ of the closure of $\fV$, where every $X \in \widetilde{\ol{\fV}}$ corresponds to a point evaluation $\Phi_X \colon f \mapsto  \widetilde{f}(X)$. 
For $H^\infty(\fV)$, we concentrate on the
weak-$*$ continuous finite-dimensional
representations, which are shown to correspond
bijectively to point evaluations $\Phi_X$
at points $X\in\widetilde{\fV}$.
We also introduce the pseudo-hyperbolic metric on the varieties which will play a role in the classification scheme later on. 

We use the representation spaces and their metric structure to obtain the basic classification theorems, Theorem \ref{thm:bounded-cb-classification} and Corollary \ref{cor:Hinf_class}, which say that $A(\fV)$ and $A(\fW)$ (or $H^\infty(\fV)$ and $H^\infty(\fW)$) are boundedly isomorphic if and only if there exists a bi-Lipschitz NC biholomorphism between the similarity envelopes $\widetilde{\fV}$ and $\widetilde{\fW}$. 
In the case of uniformly continuous function algebras the basic theorem requires also that this mapping is uniformly continuous, however, this condition of uniform continuity is later eliminated in the important Theorem \ref{thm:linearisation_maps} which gives a three way equivalence, saying that for homogeneous $\fV$ and $\fW$ the following are equivalent: 
\begin{enumerate}
\item $A(\fV)$ and $A(\fW)$ are boundedly isomorphic. 
\item There exists a bi-Lipschitz NC biholomorphism between the similarity envelopes $\widetilde{\fV}$ and $\widetilde{\fW}$. 
\item There exists a linear map $E \to F$ which induces a bi-Lipschitz bijection between the similarity envelopes $\widetilde{\fV}$ and $\widetilde{\fW}$. 
\end{enumerate}
Similar results for both the bounded and uniformly continuous function algebras are obtained also in the completely bounded category.

The classification results for homogeneous
varieties depend on
Theorems \ref{thm:classification_linear}
and \ref{thm:classification_linear_weak_star},
which establish the equivalence between
the existence of a bounded isomorphism
and that of a suitable bi-Lipschitz
NC linear bijection between the corresponding
similarity envelopes.
The linearization results rely on a new ``two-circle trick," which replaces the classical ``disc trick" for homogenizing isomorphisms; see Lemma \ref{lem:two_circ}. 
This technique is also applied to obtain a rigidity result, Theorem \ref{thm:ball_alg_iso}, saying that if there exists a completely bounded isomorphism $\alpha \colon A(\bB_E) \to A(\bB_F)$, then the operator spaces $E$ and $F$ are completely isometric. 
In Section \ref{sec:comm_balls} we show that this rigidity does not hold for the commutative slices. 
In Section \ref{sec:examples} we present a series of examples showing how our results are sharp in some sense.

The results of this paper can be seen either as generalizations to the case of arbitrary operator space balls of classification results for algebras of bounded NC functions on subvarieties of the row ball \cite{SSS1,SSS2}, or as a (completely) bounded analogue of the completely isometric classification results in the case of arbitrary operator spaces \cite{SampatShalit25,sampat2025weak}. 
It is significant and interesting to note that there are surprises: many of the results from those settings fail to generalize or to extend to our setting; see in particular Example \ref{ex:non_singleton_fiber} for non-singleton fibers of representations, or Section \ref{sec:examples} for more examples which delineate the boundaries of our positive results. 
On the other hand, new techniques are required to establish the generalizations of the results that do generalize (see for example the use of the spectral radius in the proof of Proposition \ref{prop:homo_ind_NCmap}, or the linearization results of Section \ref{sec:linear}). 
Moreover, one of the most striking theorems that we obtain --- the rigidity of NC ball algebras (Theorem \ref{thm:ball_alg_iso}) --- uncovers a phenomenon that cannot even be formulated without considering the generality of allowing general operator balls. 

\subsection*{AI disclosure.} Various AI tools were used continuously throughout the work for various routine tasks such literature search or as a thesaurus, or to expedite testing ideas and examples. 
It was also used to proofread the paper. 
We wish to single out a few outstanding points where the AI tool made the decisive contribution. First, the following examples were found by the AI system: the example following Proposition \ref{prop:rep_AV}, Example \ref{ex:non_singleton_fiber}, Example \ref{ex:cntr_exp_BLip} (ChatGPT 6 Astra), and Example \ref{ex:cb_not_ciso} (ChatGPT 5.6 Sol).
Significantly, the proof of Lemma \ref{lem:two_circ} was discovered by AI (ChatGPT 5.6 Sol) after we asked it if such a statement, which we have been trying to prove for a long time, is true. The strengthened statement of Theorem \ref{thm:linearisation_maps} and the key ideas underlying its proof were developed primarily with the assistance of ChatGPT 5.6 Sol.
The authors of the paper have envisioned, initiated, digested, understood and written up everything; we take responsibility for the correctness and we also take credit for the discoveries. 

\section{Preliminaries}

\subsection{Background on NC function theory}

There are many approaches to noncommutative (NC) analytic functions, for example Taylor \cite{Taylor-frame, Taylor-ncfunc}, Voiculescu \cite{Voiculescu-quest1, Voiculescu-quest2} and Muhly--Solel \cite{MS04,MS08,MS11}. 
We work in the framework that crystallized about a decade ago with the works of Kaliuzhnyi-Verbovetskyi and Vinnikov \cite{KVV}, Helton--Klep--McCullough \cite{HKM11}, Agler--McCarthy \cite{agler2018global}, Ball--Bolotnikov \cite{ball2021noncommutative}, Popescu \cite{Popescu-funcI,Popescu-nc_dom}  and others \cite{agler2020operator,AHKM,BMV-int}, but drawing most directly on \cite{SampatShalit25,sampat2025weak} (see also \cite{kumarsampat2025zeros,sampat2026cyclicity}).

\subsubsection{NC sets and NC domains and operator balls}

Throughout this paper $E$ and $F$ will denote finite-dimensional operator spaces. 
Let $M_n$ denote the $n \times n$ matrices over $\bC$. 
For a finite-dimensional operator space $E$, we let $M_n(E) = M_n \otimes E$ and define $\bM(E) = \cup_{n=1}^\infty M_n(E)$. We remind the reader that an operator space $E$ is just a subspace of a C*-algebra with the inherited norm. 
This implies that $M_n(E)$ too is an operator space.
We sometimes identify $E = \bC^d$ equipped with a particular operator space structure (i.e. a family of matrix norms), in which $M_n(\bC^d) = M_n^d$ is the set of all $d$-tuples of scalar matrices, and we then write $\bM^d$ for $\bM(\bC^d) = \cup_{n=1}^\infty M_n^d$. 
A subset $\Omega \subseteq \bM(E)$ is said to be an {\em NC set} if it is closed under direct sums. 
We write $\Omega_n$ or $\Omega(n)$ for $\Omega \cap M_n(E)$, so $\Omega = \cup_{n=1}^\infty \Omega_n$. 

An {\em NC domain} is an open and connected NC set. 
Every operator space $E$ gives rise to the {\em NC operator ball} $\bB_E$
\begin{equation}\label{eq:op_ball}
\bB_E = \cup_{n=1}^\infty \{X \in M_n(E) : \|X\|_n<1\}.
\end{equation}
When $E$ is an operator structure on $\bC^d$ determined by an injective linear embedding $\bC^d \hookrightarrow B(\cH)$ given by a map $e_j \mapsto Q_j \in B(\cH)$, we identify 
\begin{equation}\label{eq:op_ball_tuples}
\bB_E = \cup_{n=1}^\infty \{X \in M_n^d : \|X\|_Q := \|\sum_j X_j \otimes Q_j \| <1 \}.
\end{equation} 
It is useful to have a couple of concrete examples at hand. 
The {\em row ball} $\fB_d$, which is defined to be $\bB_E$ for $E=\bC^d$ with the row operator space structure, is given by
\begin{equation}\label{eq:Ball}
\fB_d := \bB_{(\bC^d)_{row}} = \{ X \in \bM^d : \|X\|_{row} := \|\sum X_j X_j^*\|^{1/2} < 1\} .
\end{equation}
When $E = \min(\ell^\infty_d)$ is the minimal operator space over $\ell^\infty_d$ we get the {\em NC polydisc}
\[
\fD_d := \bB_{\min(\ell^\infty_d)} = \{X \in \bM^d : \|X\|_\infty := \max_{1 \leq i \leq d} \|X_i\| < 1\}. 
\]

\subsubsection{NC analytic functions}

A function $f$ from an NC set $\Omega \subseteq \bM(E)$ to $\mathbb{M}(F)$ is said to be an {\em NC function} if: (i) $f$ is graded: $X \in \Omega_n \Rightarrow f(X) \in M_n(F)$; (ii) $f$ respects direct sums: $f(X \oplus Y) = f(X) \oplus f(Y)$; and (iii) $f$ respects similarities: if $X \in \Omega_n$, $S \in GL_n$, and if $S^{-1} X S \in \Omega_n$, then $f(S^{-1} X S) = S^{-1} f(X) S$. 
For example, {\em free polynomials} are NC functions on $\bM^d$: if $p(Z) = \sum_w c_w Z^w$ a free polynomial, then one evaluates $p$ on $A \in \bM^d$ as follows
\[
p(A) = \sum_w c_w A^w, 
\]
where for a word $w$ in $\{1,2,\ldots, d\}^n$ the monomial $Z^w = Z_{w_1} \cdots Z_{w_n}$ evaluated at $A^w$ is given by $A_{w_1} \cdots A_{w_n}$. 
Below we shall always let $Z_1,\ldots,Z_d$ denote the natural coordinate functions on $\bM^d$.  
A more general class of NC functions is provided by {\em NC rational functions} \cite{helton2018applications,JMS-rat,KVV-rat,mai2018free,Volcic}. 
Remarkably, every bounded NC function is {\em analytic} and has a Taylor series at every point; in fact, for a bounded NC function in an operator ball, the series around the origin converges in the entire ball \cite{KVV}. 

The {\em similarity envelope} $\widetilde{\Omega}$ of an NC set $\Omega$ consists of all tuples $S^{-1} X S := (S^{-1}X_1S, \ldots, S^{-1}X_d S)$ where $X \in \Omega$, that is
\begin{equation}\label{eq:sim_env}
\widetilde{\Omega} = \bigcup_{n=1}^\infty\left\{S^{-1} X S : X \in \Omega_n, \,\, S \in \GL_n \right\}.
\end{equation}
Every NC function $f$ defined on an NC set $\Omega$ always extends to an NC function $\widetilde{f}$ on the similarity envelope $\widetilde{\Omega}$ of $\Omega$, given by the formula $\widetilde{f}(S^{-1}XS) = S^{-1}f(X)S$ (see the appendix by S. Biswas in \cite{KVV}). 

\subsubsection{NC varieties and the algebras of interest}

Let $\Omega$ be an NC domain. 
An {\em NC analytic variety $\fV \subseteq \Omega$} (or simply a {\em variety}, for short) is defined to be the joint zero set of a set of bounded NC functions on $\Omega$.   
A variety is {\em homogeneous} if it the joint zero set of a homogeneous ideal in $\bC\langle Z_1, \ldots, Z_d\rangle$, or, equivalently, a homogeneous ideal in $A(\bB_E)$ or $H^\infty(\bB_E)$. A variety $\fV$ is homogeneous if and only if $tX \in \fV$ for all $t \in \bD$ whenever $X \in \fV$ (see \cite{SampatShalit25}). 

We define $H^\infty(\fV)$ to be the algebra of bounded scalar-valued NC functions on $\fV$ equipped with the supremum norm
\[
\|f\| = \|f\|_\infty := \sup_{X \in \fV}\|f(X)\|. 
\] 
We define $A(\fV)$ to be the closure of polynomials in $H^\infty(\fV)$. In particular, when $\Omega = \fV = \bB_E$ we have the algebras $A(\bB_E)$ and $H^\infty(\bB_E)$ at the center of the stage. 
If $\fV$ is a homogeneous set, then $A(\fV)$ equals the set of NC functions that are uniformly continuous on ${\fV}$ (see \cite[Proposition 4.12]{SampatShalit25}); equivalently, these are the bounded NC functions that extend uniformly continuously to $\ol{\fV}$.
For notational convenience, we let $\cA(\fV)$ denote {\em either} $A(\fV)$ {\em or} $H^\infty(\fV)$, and refer to both as {\em NC function algebras}. 
It is easy to verify that $\cA(\fV)$, equipped with the supremum norm and with the pointwise operations of addition and multiplication, satisfies the Blecher--Ruan--Sinclair axioms and so is an operator algebra \cite{BlecherRuanSinclair1990}.

\subsubsection{Homogeneous decomposition}

Every $f \in H^\infty(\bB_E)$ can be written formally as a series of homogeneous polynomials
\begin{equation}\label{eq:homog_exp}
f \sim \sum_{n=0}^\infty f_n, 
\end{equation}
where the maps $P_n:H^\infty(\bB_E)\longrightarrow H^\infty(\bB_E), \, P_n f=f_n,$
are given explicitly by 
\begin{equation}\label{eq:homog_proj}
P_n f(X) = f_n(X) = \frac{1}{2\pi} \int_0^{2\pi} f(e^{it}X)e^{-int} dt \quad , \quad \textrm{ for } X \in \bB_E,
\end{equation}
and are completely contractive for all $n \in \bN$.
The series \ref{eq:homog_exp} converges pointwise in $\bB_E$, Ces\`{a}ro converges in the bounded pointwise topology in $H^\infty(\bB_E)$, and Ces\`{a}ro converges in norm if $f \in A(\bB_E)$. 
In other words, the Ces\`{a}ro sums
\begin{equation}\label{eq:Cesaro}
    \Sigma_N(f) = \sum_{n=0}^{N-1} \left( 1 - \frac{n}{N}\right) f_n \xrightarrow{N \to \infty} f
\end{equation}
in the corresponding topology. 
Moreover, the decomposition \eqref{eq:homog_exp}  and the convergence \eqref{eq:Cesaro} also hold for $f \in H^\infty(\fV)$ or $f\in A(\fV)$ when $\fV$ is a homogeneous variety; see \cite[Proposition 3.5]{SampatShalit25} for details.

\subsubsection{Restriction and quotient}\label{subsubsection:quotient}

Let $\fV\subseteq \bB_E$ be an NC variety. Consider the restriction maps
$H^\infty(\bB_E)\longrightarrow H^\infty(\fV)$
and
$A(\bB_E)\longrightarrow A(\fV)$
and their kernels
\begin{equation}\label{eq:JequalIV}
J_\fV
=
\left\{
g\in H^\infty(\bB_E):
g(Y)=0 \textrm{ for all } Y\in\fV
\right\},
\end{equation}
and
$I_\fV
=
J_\fV\cap A(\bB_E).$
By \cite[Theorem 4.8]{SampatShalit25}
\[
H^\infty(\fV)
\cong
H^\infty(\bB_E)/J_\fV 
\]
completely isometrically. If one assumes that $\fV$ is homogeneous, then also 
\[
A(\fV)
\cong
A(\bB_E)/I_\fV
\]
completely isometrically. 
Moreover, by ``the trivial Nullstellensatz" (\cite[Proposition 4.6]{SampatShalit25}), the variety $\fV$ is recovered as the common zero set of $J_\fV$:
\begin{equation}\label{eq:VequalZJ}
\fV
=
\left\{
Y\in\bB_E:
g(Y)=0 \textrm{ for all } g\in J_\fV
\right\}.
\end{equation}

\subsubsection{Dual algebra structure}
The algebras $H^\infty(\fV)$ are endowed with a natural notion of pointwise convergence. 
In \cite[Section 3]{sampat2025weak} it was shown that $H^\infty(\fV)$ has a natural predual so that convergence in the weak-$*$ topology in $H^\infty(\fV)$ is the same as pointwise convergence on bounded sets. 
Thus, the Ces\`{a}ro sums \eqref{eq:Cesaro} of a function $f \in H^\infty(\bB_E)$ converge in the weak-$*$ topology to $f$.

\subsection{Background on the spectral radius}

We recall the spectral radius for matrices over an operator space introduced in \cite{shalsham2025spectral}. 
Let $E$ be a finite-dimensional operator space. For
$X\in B(\mathcal K)\otimes_{\min}E$, the $E$-spectral radius is defined by
\[
\rho_E(X)
:=
\rho(\iota(X))
=
\lim_{m\to\infty}\|\iota(X)^m\|^{1/m},
\]
where
$\iota:
B(\mathcal K)\otimes_{\min}E
\longrightarrow
B(\mathcal K)\otimes_{\min}C^*_{\max}(E)$
is the canonical embedding. 
If $E$ is an operator space structure on $\mathbb C^d$ and
$X=(X_1,\ldots,X_d)\in M_n^d$, we write $\rho_E(X)$ for the spectral
radius of the corresponding element of $M_n(E)$. The fundamental
characterization is
$\rho_E(X)<1$ if and only if
$X\in\widetilde{\bB}_E,$
or equivalently, there exists $S\in GL_n(\mathbb C)$ such that
$\|S^{-1}XS\|_E<1$ (see also \cite{SchShaSha} for a proof of this fact in the case of infinite-dimensional operator spaces). 
Moreover,
\[
\rho_E(X)
=
\inf_{S\in GL_n(\mathbb C)}
\|S^{-1}XS\|_E.
\]

If $X=(X_1,\ldots,X_d)$ is a commuting tuple, then
$\rho_E(X)
=
\max_{\lambda\in\sigma(X)}\|\lambda\|_E,$
where $\sigma(X)\subseteq\mathbb C^d$ is the joint spectrum \cite{SchShaSha}. 
Hence,
$\rho_E(X)<1$ if and only if
$\sigma(X)\subseteq\bB_E(1)$. 
Thus, for commuting tuples, $\rho_E(X)<1$ depends only on the first level
of the operator ball.
On the other hand, if two operator space structures $E$ and $F$ on $\mathbb C^d$ satisfy $\rho_E(X)=\rho_F(X)$ for every $d$-tuple $X$ of bounded operators on every Hilbert space, then the identity map $E\to F$ is completely isometric \cite{tornes2026} (see also \cite{SchShaSha}). 

We shall need the following Schwarz-type lemma for the spectral radius. 
If $G$ is an NC function in an NC neighborhood of $0$, we let $\Delta G(0,0)$ denote the difference differential operator of $G$ at $(0,0)$; see \cite[Section 2.2]{KVV}. 
The property of $\Delta G(0,0)$ that we shall use below is that $X \mapsto \Delta G(0,0)[X]$ is a linear map that satisfies 
\[
\Delta G(0,0)[X] = \lim_{z \mapsto 0} \frac{1}{z} \left(G(zX) - G(0)\right)
\]

\begin{lemma}\label{lem:spectral_less}
Let $E$ and $F$ be two finite-dimensional operator space structures. 
If $G \colon \widetilde{\bB}_E \to \widetilde{\bB}_F$ is NC holomorphic and satisfies $G(0) = 0$, then 
\[
\rho_F(\Delta G(0,0)[X]) \leq \rho_E(X), 
\]
for all $X \in \widetilde{\bB}_E$. 
\end{lemma}
\begin{proof}
By a standard scaling argument, it suffices to show that $\rho_E(X) < 1$ implies that $\rho_F(\Delta G(0,0)[X]) \leq 1$. 
Fix $X \in \widetilde{\bB}_E$; say $X$ is in level $n$, i.e. $X \in M_n(E)$. 
Define a function $g$ from a neighborhood of $\ol{\bD}$ into $M_n(F)$ by 
\[
g(z) = 
\begin{cases}
\frac{1}{z} G(zX) \, & \, z \neq 0\\  
\Delta G(0,0)[X] \, & \, z = 0.
\end{cases}
\]
Then $g$ is holomorphic, and therefore, by \cite[Lemma 2.5]{shalsham2025spectral},
\[
u(z) = \rho_F(g(z))
\]
is subharmonic. 
Since $G$ maps $\widetilde{\bB}_E$ into $\widetilde{\bB}_F$, for every $z \in \partial \bD$ we have $u(z) < 1$. 
By the maximum principle, 
\[
\rho_F(\Delta G(0,0)[X]) = \rho_F(g(0)) = u(0) \leq 1, 
\]
as required. 
\end{proof}

\section{The space of finite-dimensional representations}\label{sec:fin_reps}

For a unital operator algebra $\mathcal A$, let
$\operatorname{Rep}_{\mathrm{fin}}(\mathcal A)$ denote the collection
of all bounded finite-dimensional representations of $\mathcal A$; that is,
\[
\operatorname{Rep}_{\mathrm{fin}}(\mathcal A)
:=
\bigcup_{n\geq 1}
\left\{
\Phi:\mathcal A\longrightarrow M_n(\mathbb C):
\Phi \text{ is a unital bounded algebra homomorphism}
\right\}.
\]
For each $n\geq 1$, we write
$\operatorname{Rep}_n(\mathcal A)
:=
\left\{
\Phi\in \operatorname{Rep}_{\mathrm{fin}}(\mathcal A):
\Phi:\mathcal A\longrightarrow M_n(\mathbb C)
\right\}.$
Thus,
\[
\operatorname{Rep}_{\mathrm{fin}}(\mathcal A)
=
\bigcup_{n\geq 1}\operatorname{Rep}_n(\mathcal A).
\]
Throughout this paper, and in particular in this section, $E$ will denote a finite-dimensional operator space given by a concrete operator space structure on $\bC^d$ and $Z_1, \ldots, Z_d$ will denote the coordinate functions on $\bB_E$. 

\subsection{Representations of algebras of uniformly continuous NC functions}

Let $\fV \subseteq \bB_E$ be a homogeneous variety. 
Then, for every $n$ and every $X\in\ol{\fV}(n)$, evaluation at $X$ gives rise to a representation $\Phi_X:A(\fV)\longrightarrow M_n(\mathbb C)$ given by 
\[
\Phi_X ( f)=  f(X),
\]
where we use the continuous extension of $f\in A(\fV)$ to
$\ol{\fV}$.
Since every NC function $f$ defined on an NC set $\Omega$ always extends to an NC function $\widetilde{f}$ on the similarity envelope $\widetilde{\Omega}$, every
$X\in\widetilde{\ol{\fV}}(n)$ gives rise to a bounded evaluation representation
$\Phi_X:A(\fV)\longrightarrow M_n(\mathbb C)$ given by $\Phi_X(f)=\widetilde{f}(X)$.

\begin{proposition}\label{prop:rep_AV}
    Let $\fV \subseteq \bB_E$ be a homogeneous variety. There is a bijective correspondence between $\widetilde{\ol{\fV}}$ and bounded finite-dimensional representations of $A(\fV)$ given by 
    \[
    X \mapsto \Phi_X
    \]
    with inverse given by 
    \begin{equation}\label{eq:def_pi}
    \Phi \mapsto \pi(\Phi) := (\Phi(Z_1), \ldots, \Phi(Z_d)), 
    \end{equation}
    where $Z_1,\ldots,Z_d$ are the coordinate functions.
\end{proposition}

\begin{proof}
Define $\pi$ by \eqref{eq:def_pi}. 
We first show that $\pi$ maps into $\widetilde{\ol{\fV}}$.
Indeed, if 
$\Phi\in\operatorname{Rep}_{\mathrm{fin}}(A(\fV))$, then, for some $n\in\mathbb N$,
$\Phi:A(\fV)\longrightarrow M_n(\mathbb C)$
is a unital bounded algebra homomorphism.
Smith's theorem
implies that $\Phi$ is completely bounded.
By Paulsen's similarity theorem
\cite{PaulsenBook}, there exists an invertible
$S\in\GL_n(\mathbb C)$ such that
$\rho(f)=S^{-1}\Phi(f)S$
is a unital completely contractive representation of $A(\fV)$.
By \cite[Proposition 5.3]{SampatShalit25}, there exists a unique
$X\in\ol{\fV}(n)$ such that
$\rho=\Phi_X$. Since every $f\in A(\fV)$ extends continuously to
$\ol{\fV}$, we use the same notation $f$ for this extension. Set
$Y:=SXS^{-1}\in\widetilde{\ol{\fV}}(n).$
Consequently, 
\[
\Phi(f)
=
S\rho(f)S^{-1}
=
Sf(X)S^{-1}
=
\widetilde f(Y)
=
\Phi_Y(f),
\qquad
f\in A(\fV).
\]
Therefore,
\[
\pi(\Phi)
=
\bigl(\Phi_Y(Z_1),\ldots,\Phi_Y(Z_d)\bigr)
=
\bigl(SX_1S^{-1},\ldots,SX_dS^{-1}\bigr)
=
Y.
\]
Thus every
$\Phi\in\operatorname{Rep}_{\mathrm{fin}}(A(\fV))$
corresponds to a point of
$\widetilde{\ol{\fV}}$.
To show that $\pi$ is surjective, let
$Y\in\widetilde{\ol{\fV}}$. 
As observed above, we get $\Phi_Y\in\operatorname{Rep}_{\mathrm{fin}}(A(\fV))$ such that $\Phi_Y(Z_j)=Y_j,
\, j=1,\ldots,d.$ Therefore
$\pi(\Phi_Y)
=
Y$. 
Finally, $\pi$ is injective since $A(\fV)$ is generated by the coordinate functions. 
\end{proof}

Although we shall not use topological considerations on the space of representations, it is natural to wonder whether the above correspondence is continuous. 
If we equip $\operatorname{Rep}_n(A(\fV))$ with the topology of pointwise convergence (i.e. $\Phi_\alpha \to \Phi$ if and only if $\Phi_\alpha(f) \to \Phi(f)$ for all $f \in A(\fV)$) and $\widetilde{\ol{\fV}}(n)$ with the topology inherited from $M_n^d$, then it is clear from \eqref{eq:def_pi} that $\pi$ is continuous. 
The inverse $\pi^{-1}$ is continuous when restricted to $\widetilde{\fV}(n)$. 
Indeed, if $X_k \to X \in \widetilde{\fV}$ then $S^{-1}XS \in \fV$ for some invertible $S$ and then $S^{-1} X_k S \in \fV$ too for sufficiently large $k$ and $S^{-1}X_k S \to S^{-1} X S$. 
Thus, for all $f \in A(\fV)$
\[
\Phi_{X_k}(f) = Sf(S^{-1}X_k S)S^{-1} \to Sf(S^{-1}X S)S^{-1} = \Phi_X(f), 
\]
and so $X \mapsto \Phi_X$ is continuous on $\widetilde{\fV}$. 
However, $X \mapsto \Phi_X$ need not be continuous on $\widetilde{\ol{\fV}}(n)$. 
To see this, let $\fV = \fB_1 = \bB_{\bC}$ be the one-variable ball of strict contractions, and consider 
\[
X_k = \begin{pmatrix} 1-\frac{1}{k^2} & \frac{1}{k} \\ 0 & 1-\frac{1}{k^2} \end{pmatrix} \in \widetilde{\fV}(2).
\]
Let $f \in A(\fV) = A(\bD)$ be given by $f(z) = \sqrt{1-z}$ where the square root is given by the principal branch. 
Then 
\[X_k \xrightarrow{k \to \infty} I_2 \in \widetilde{\ol{\fV}}(2)\] 
but 
\[
\Phi_{X_k}(f) = f(X_k) 
= \begin{pmatrix} \frac{1}{k} & -\frac{1}{2} \\ 0 & \frac{1}{k} \end{pmatrix} \xrightarrow{k\to \infty} \begin{pmatrix} 0 & -\frac{1}{2} \\ 0 & 0 \end{pmatrix} \neq 0 = f(I_2). 
\]
This shows that the map $X \mapsto \Phi_X$ is not continuous at the point $X = I_2$. 

\subsection{Representations of algebras of bounded NC functions}

When dealing with the full algebra of bounded NC functions $H^\infty(\fV)$ we need not restrict attention to homogeneous varieties. The representation spaces of $H^\infty$ algebras are significantly wilder. In that case, we do not have an evaluation representation for every point in $\ol{\fV}$. 
However, we clearly have an evaluation $\Phi_X$ for every $X \in \fV$ and therefore also for every $X \in \widetilde{\fV}$. 

\begin{theorem}\label{thm:weakstar_rho_lessthan_1}
Let $\fV\subseteq\bB_E$ be a variety.
There is a continuous map 
\[
\pi\colon \operatorname{Rep}_{\operatorname{fin}}(H^\infty(\fV)) \to \widetilde{\ol{\bB_E}}
\]
given by $\Phi \mapsto \pi(\Phi) := (\Phi(Z_1), \ldots, \Phi(Z_d))$, which restricts to a bijection between the weak-$*$-continuous representations and $\widetilde{\fV}$.
If $\fV$ is homogeneous, then $\pi$ is a surjection 
\[
\pi\colon \operatorname{Rep}_{\operatorname{fin}}(H^\infty(\fV)) \to \widetilde{\ol{\fV}} .
\]
\end{theorem}
\begin{proof}
The main content of the theorem is the bijective correspondence between the weak-$*$ representations in $\operatorname{Rep}_{\operatorname{fin}}(H^\infty(\fV))$ and $\widetilde{\fV}$. We dispence with the easier claims first. 

The fact that $\pi$ is continuous and maps $\operatorname{Rep}_{\operatorname{fin}}(H^\infty(\fV))$ into $\widetilde{\ol{\bB_E}}$ follow from \cite[Proposition 5.1]{SampatShalit25} by an argument analogous to that used in Proposition \ref{prop:rep_AV}. 

Suppose that $\fV$ is homogeneous. 
Every $\Phi \in \operatorname{Rep}_{\operatorname{fin}}(H^\infty(\fV))$ restricts to a representation of $A(\fV)$, so by Proposition \ref{prop:rep_AV} $X = \pi(\Phi) \in \widetilde{\ol{\fV}}$. 
Conversely, if $X  \in \widetilde{\ol{\fV}}$, then $rX \in \widetilde{\fV}$ for all $r \in (0,1)$ and gives rise to the point evaluation representation $\Phi_{rX} \in \operatorname{Rep}_{\operatorname{fin}}(H^\infty(\fV))$. 
It follows that $rX$ is in the range of $\pi$ and we conclude that $\widetilde{\fV}$ is in the range of $\pi$. 
A compactness argument shows that $\pi$ surjects onto $\widetilde{\ol{\fV}}$, as stated. 

We now show that the map $\pi$ identifies the weak-$*$ continuous finite-dimensional representations of $H^\infty(\fV)$ with the points of $\widetilde{\fV}$.
First, every $X\in\widetilde{\fV}$ gives rise to a weak-$*$ continuous point evaluation representation
$\Phi_X$, and
$\pi(\Phi_X)=X$. 
Conversely, let
$\Phi\in\operatorname{Rep}_{\operatorname{fin}}(H^\infty(\fV))$
be weak-$*$ continuous and set
$X=\pi(\Phi)$.
Since $\Phi$ lifts to a weak-$*$ continuous representation of $H^\infty(\bB_E)$ and every $f \in H^\infty(\bB_E)$ is the bounded pointwise limit of polynomials as in \eqref{eq:Cesaro}, once it is shown that $X\in\widetilde{\fV}$ it would follow that
$\Phi(f)=\Phi_X(f)
$ for every $f\in H^\infty(\fV)$.
It therefore remains to prove the essential point that
$X\in\widetilde{\fV}$.

We will show below that $\rho_E(X)<1$. 
Let us show how this implies that $X \in \widetilde{\fV}$. 
Indeed, suppose without loss of generality that $X \in \bB_E$. Since the lift of $\Phi$ to $H^\infty(\bB_E)$ has to annihilate the ideal $J_\fV = \{g \in H^\infty(\bB_E) : g\big|_\fV \equiv 0\}$, we have $g(X) = 0$ for all $g \in J_\fV$, which means that $X \in \fV$; see Equations \eqref{eq:JequalIV} and \eqref{eq:VequalZJ}.

We now show that $\rho_E(X) < 1$. 
Since $\Phi$ lifts to a weak-$*$ continuous representation of $H^\infty(\bB_E)$, we may as well assume that $\fV = \bB_E$. 
Therefore, assume for the sake of obtaining a contradiction that $\Phi$ is a weak-$*$ continuous representation of $H^\infty(\bB_E)$ such that $X = \pi(\Phi)$ and that $\rho_E(X) = 1$. 
Without loss of generality, we may assume $X$ is irreducible (see \cite[Remark 2.8]{shalsham2025spectral} and then by \cite[Corollary 2.12]{shalsham2025spectral} we may further assume that $\|X\| = 1$. 

Consider the linear pencil $\sum_j X_j \otimes W_j$ as an NC function on the unit ball $\bB_{E^*}$ of the operator space dual $E^*$. 
Since the coordinate functions $W_1, \ldots, W_d$ can be identified with the generators $Q_1, \ldots, Q_d$ of $E$, we have that $\rho(\sum X_j \otimes W_j) = \rho(\sum X_j \otimes Q_j) = \rho_E(X) = 1$, by assumption. Note that $\left\| \sum X_j \otimes W_j \right\| = 1$. 
By \cite[Corollary 3.3]{shalsham2025spectral}, there exists some $n \in \bN$ and $Y \in \ol{\bB_{E^*}}(n)$ such that 
\begin{equation}\label{eq:1inspec}
1 \in \sigma\left(\sum_{j=1}^d X_j \otimes Y_j \right) .
\end{equation}
Now consider the $M_n$ valued NC function $f \colon \bB_E \to \bM(M_n)$ by 
\[
f(Z) = \exp\left\{i \operatorname{Log}\left[\left(1 + \sum_j Z_j \otimes Y_j \right)\left(1 - \sum_j Z_j \otimes Y_j\right)^{-1} \right] \right\} = \phi\left( \sum Z_j \otimes Y_j\right), 
\] 
where 
\[
\phi(t) = \exp\left(i \operatorname{Log} \left(\frac{1+t}{1-t} \right) \right)
\]
is a bounded analytic function on $\bD$ that has no radial limit at $1$. 
By the fact that $\left\| \sum Z_j \otimes Y_j \right\| \leq 1$ together with von Neumann's inequality for matrix valued functions, $f \in M_n(H^\infty(\bB_E))$ and $\|f\| \leq e^{\pi/2}$. 
Note that we are not claiming that $f$ may be evaluated at $X$. 

For $r \in (0,1)$ we define $f_r(Z) = f(rZ)$. 
It is easy to see that $f_r \xrightarrow{r \to 1^-} f$ bounded-pointwise, equivalently, in boundedly in the weak-$*$ topology. 
Thus, 
\[
f_r(X) = \Phi^{(n)}(f_r) \longrightarrow \Phi^{(n)}(f) \textrm{ as } r \to 1^{-}
\]
for all $f$, where the equality follows from $\pi(\Phi) = X$ and the uniform convergence of the Taylor series of $f_r$. 
On the other hand, Equation \eqref{eq:1inspec} shows that that at least one of the entries of the matrix $f_r(X)$ does not converge as $r \to 1^{-}$. 

This contradiction leads to the conclusion that we cannot have both $\Phi$ weak-$*$ continuous and $\rho_E(\pi(\Phi)) = 1$, and so we conclude that $\rho_E(X)<1$ for every $X = \pi(\Phi)$ which is the projection of a weak-$*$ continuous representation $\Phi$. 
\end{proof}

\begin{remark}
Let us present an alternative way to prove that if $\Phi \in \operatorname{Rep}_{\operatorname{fin}}(H^\infty(\bB_E))$ is weak-$*$ continuous then $\rho_E(\pi(\Phi)) < 1$. 
This proof is based on a direct examination of the definition of spectral radius $\rho_E(X)$. 

Let $\Phi$ be such a representation and put $X = \pi(\Phi) \in M_k^d$. 
Following the notation in the proof of \cite[Theorem 2.7]{shalsham2025spectral}, we define for every $n$ the space of $n$-homogeneous polynomials
\[
E_n = \spn\{Z^w : |w| = n\}. 
\]
We define $\pi_{n,X} \colon E_n \to M_k$ by point evaluation $\pi_{n,X}(p) = p(X)$. 
By the definition of $\rho_E$, we have that 
\[
\rho_E(X) = \lim_{n\to \infty} \|\pi_{n,X}\|_{cb}^{1/n} \,\, , 
\]
see \cite{shalsham2025spectral} for details. 
We shall show that $\lim_{n\to \infty} \|\pi_{n,X}\|_{cb} = 0$. Since $\|\pi_{n,X}\|_{cb}$ is a sub-multiplicative sequence, it would then follow by Fekete's lemma that 
\[
\rho_E(X) = \lim_{n\to \infty} \|\pi_{n,X}\|_{cb}^{1/n} = \inf_n \|\pi_{n,X}\|_{cb}^{1/n} < 1. 
\]
Since $\pi_{n,X}$ maps into $M_k$, it follows that $\|\pi_{n,X}\|_{cb} = \|\pi_{n,X}\|_k$. 
Thus, for every $n$, there exists $P_n \in M_k(E_n)$ with $\|P_n\| = 1$ such that $\|\pi_{n,X}\|_{cb} = \|\pi^{(k)}_{n,X}(P_n)\| = \|P_n(X)\|$. 
But it is easy to see that $P_n$, as a sequence of norm $1$ homogeneous polynomial of degree $n$, converges to $0$ in the bounded pointwise topology. 
Since $\Phi$ is weak-$*$ continuous, we have
\[
\|\pi_{n,X}\|_{cb} = \|P_n(X)\| = \|\Phi^{(k)}(P_n)\| \to 0, 
\]
as required. 
\end{remark}

\subsection{Non-singleton fibers --- an example}

Experience with classical function algebras, as well as with algebras of bounded NC functions on the row ball, suggests a kind of converse to the main result of the previous subsection. Namely, one might think that if a representation $\Phi \in \operatorname{Rep}_{\operatorname{fin}}(H^\infty(\bB_E))$ is fibered over the similarity envelope of $\bB_E$ (that is, $\rho_E(\pi(\Phi)) < 1$), then $\Phi$ must be weak-$*$ continuous and therefore equal to $\Phi_X$ for some $X \in \widetilde{\bB}_E$. 
For some time it was believed that this is the case. 
Efforts to prove this result in the general case resulted in the identification of function theoretic obstructions to certain proof strategies (see \cite[Section 5]{sampat2025weak}) but not to a proof. 
We now resolve this problem by showing that there do exist representations fibered over the interior of the ball which are not weak-$*$ continuous. 

\begin{lemma} \label{lem:coeff_diff}
    Let $J_n$ be the nilpotent Jordan block of size $(n+1) \times (n+1)$. Write $\{v_0,\ldots,v_n\}$ for the standard basis of $\C^{n+1}$. Let $J'_n = \left( \begin{smallmatrix} -1 & 0 \\ 0 & I_n \end{smallmatrix} \right) J_n$. Then, for every homogeneous polynomial $f(Z_1, Z_2) = \sum_{|\alpha| = m} a_{\alpha} Z^{\alpha} \in \C \langle Z_1, Z_2 \rangle$ with $\deg f = m \geq 1$,
    \begin{itemize}
        \item[(i)] If $m > n$, then $f(J_n,J'_n) = 0$;

        \item[(ii)] If $1 \leq m \leq n$ and $1 \leq j < m$, then $f(J_n,J'_n) v_j = 0$;

        \item[(iii)] If $1 \leq m \leq n$ and $m < j \leq n$, then $f(J_n, J'_n)v_j = \left( \sum_{|\alpha| = m} a_{\alpha} \right) v_{j-m}$;

        \item[(iv)] If $1 \leq m \leq n$, $f(J_n,J'_n) v_m = \sum_{|\beta| = m-1} (a_{1 \beta} - a_{2 \beta}) v_0$.

    \end{itemize}
    We conclude that for every $f \in \C \langle Z_1, Z_2 \rangle$, if $f = \sum_{k=0}^m \sum_{|\alpha|=k} a_{\alpha} Z^{\alpha}$, then 
    \[
    f(J_n,J'_n)^* v_0 = \overline{f(0)} v_0 + \sum_{k=1}^n \left(\sum_{|\beta|=k-1} \overline{(a_{1\beta} - a_{2\beta})} \right) v_k.
    \]
\end{lemma}
\begin{proof}
    The last claim follows immediately from the items. Indeed, if we write $f = \sum_{k=0}^m f_k$, the homogeneous decomposition of $f$, then
    \[
    f(J_n,J'_n)^* v_0 = \sum_{k=0}^n \langle v_0, f(J_n,J'_n) v_k \rangle v_k.
    \]
    Since the point $(J_n,J'_n)$ is jointly nilpotent of order $n+1$, the claim (i) is trivial. The claims (ii) is similarly easy because of the Jordan block form. For claims (iii) and (iv), we note that it suffices to prove them for monomials. Therefore, let $f = Z^{\alpha}$, where $1 \leq |\alpha| = m \leq n$. We note that for $2 \leq i \leq n$, $J_n v_i = J'_n v_i$. Therefore, for every $m < j \leq n$, we get that $f(J_n, J'_n) v_j = J_n^m v_j = v_{j-m}$. This proves (iii). To prove (iv), we first assume that $\alpha = 1 \alpha'$. By (iii), $(J_n,J'_n)^{\alpha'} v_m = v_1$. Hence $f(J_n,J'_n) v_m = v_0$. However, if $\alpha = 2\alpha'$, we get that $f(J_n, J'_n) v_m = J_n'v_1 = - v_0$. This proves (iv).
\end{proof}

\begin{example}\label{ex:non_singleton_fiber}
   Let us consider the NC bidisc $\fD_2$. For $m \in \bN$ consider the $M_2$-valued analytic function on $\bD$ given by
   \[
   \varphi_m(\lambda) = \begin{pmatrix} 0 & \frac{\lambda}{\sqrt{m}} \\ \frac{\lambda}{\sqrt{m}} & 1 - \frac{1}{m} \end{pmatrix}.
   \]
   We define the point $Y_m(\lambda) = (\varphi_m(\lambda), \varphi_m(-\lambda))$. Note that
   \[
    1 - \varphi_m(\lambda)^* \varphi_m(\lambda) = \begin{pmatrix} 1 - \frac{|\lambda|^2}{m} & -(1 - \frac{1}{m}) \frac{\bar{\lambda}}{\sqrt{m}} \\ -(1 - \frac{1}{m}) \frac{\lambda}{\sqrt{m}} & 1 - (1 - \frac{1}{m})^2 - \frac{|\lambda|^2}{m} \end{pmatrix}
   \]
   Note that since $|\lambda| <1 $ and $m \in \bN$, $1 - \frac{|\lambda|^2}{m} > 0$. Moreover,
   \begin{multline*}
   \det(1 - \varphi_m(\lambda)^* \varphi_m(\lambda)) = \left(1 - \frac{|\lambda|^2}{m}\right)\left(1 - (1 - \frac{1}{m})^2 - \frac{|\lambda|^2}{m} \right) - (1 - \frac{1}{m})^2 \frac{|\lambda|^2}{m} = \\ \left(1 - \frac{|\lambda|^2}{m}\right)^2 - (1 - \frac{1}{m})^2 = \left(2 - \frac{|\lambda|^2 +1}{m} \right) \frac{1 - |\lambda|^2}{m} > 0.
   \end{multline*}
   Hence $Y_m \colon \bD \to \fD_2(2)$ is an analytic map. Since evaluation at a point is a homomorphism, we get a homomorphism $F_m \colon H^{\infty}(\fD_2) \to M_2(H^{\infty}(\bD))$ given by $F_m(f)(\lambda) = f(Y_m(\lambda))$. This homomorphism is clearly completely contractive.

   Now note that
   \[
   \varphi_m(-\lambda) = \begin{pmatrix} 1 & 0 \\ 0 & -1 \end{pmatrix} \varphi_m(\lambda) \begin{pmatrix} 1 & 0 \\ 0 & -1 \end{pmatrix}.
   \]
   Therefore, for every $f \in H^{\infty}(\fD_2)$,
   \[
   F_m(f)(-\lambda) = f(Y_m(-\lambda)) = \begin{pmatrix} 1 & 0 \\ 0 & -1 \end{pmatrix} F_m(f)(\lambda)\begin{pmatrix} 1 & 0 \\ 0 & -1 \end{pmatrix}.
   \]
   We conclude that the diagonal entries of $F_m(f)$ are even and the off-diagonal ones are odd. The direct sum property of NC functions tells us that $F_m(f)(0) = f(0,0) \oplus f\left(1 - \frac{1}{m}, 1 - \frac{1}{m}\right)$. 
   
   Our goal is to construct a continuous derivation on $H^{\infty}(\fD_2)$ with values in $\C$ viewed as a bimodule over $H^{\infty}(\fD_2)$ via the homomorphism of evaluation at $(0,0)$. To this end we define for some $0 < r < 1$
   \[
   d_m(f) = \frac{1}{2\pi i} \int_{|\lambda| = r} \frac{1}{\lambda^3} \langle F_m(f)(\lambda) e_1, e_1 \rangle d \lambda.
   \]
   Clearly, $d_m(f)$ is the second Taylor coefficient of the analytic function $\langle F_m(f)(\lambda) e_1, e_1 \rangle$ at the origin. In particular, $|d_m(f)| \leq \|f\|$. Therefore, $d_m$ is a continuous linear functional on $H^{\infty}(\fD_2)$. We write $F_m(f)(\lambda) = \sum_{n=0}^{\infty} A_n^{(m)}(f) \lambda^n$. The above observation implies that that for every $k \in \bN$, $A^{(m)}_{2k}$ is diagonal and $A^{(m)}_{2k-1}$ is anti-diagonal. Since the map $F_m$ is a homomorphism, we get that for every $f, g \in H^{\infty}(\fD_2)$,
   \begin{multline*}
   F_m(fg)(\lambda) = F_m(f)(\lambda)F_m(g)(\lambda) = A^{(m)}_0(f) A^{(m)}_0(g) + \lambda( A^{(m)}_0(f) A^{(m)}_1(g) + \\ A^{(m)}_1(f) A^{(m)}_0(g)) + \lambda^2 ( A^{(m)}_0(f) A^{(m)}_2(g) + A^{(m)}_1(f) A^{(m)}_1(g) + A^{(m)}_2(f) A^{(m)}_0(g)) + \cdots.
   \end{multline*}
   Therefore,
   \begin{equation} \label{eq:d_m_product}
   d_m(fg) = f(0,0) d_m(g) + d_m(f) g(0,0) + \langle A_1^{(m)}(f) e_2, e_1 \rangle \langle A_1^{(m)}(g) e_1, e_2 \rangle.
   \end{equation}
   Our candidate for our derivation is a weak* accumulation point of the $d_m$. Therefore, we need to show that the third term above goes to $0$, as $m \to \infty$. To this end, consider $f = Z^{\alpha}$ a monomial. Clearly, $F_m(f)(\lambda)$ is a matrix valued polynomial. We note that $F_m(1) = 1$, $F_m(Z_1) = \varphi_m(\lambda)$, and $F_m(Z_2) = \varphi_m(-\lambda)$. Write $f(Z_1,Z_2) = \sum_{n=0}^{\infty}f_n(Z_1,Z_2)$, where each $f_n$ is the $n$-th homogeneous component of the Taylor-Taylor expansion of $f$ around the origin. The series converges absolutely and uniformly on $r \overline{\fD_2}$ for every $0 < r < 1$.  Now for every word $\alpha$, we see that 
   \[
   F_m(Z^{\alpha})(\lambda) = \begin{pmatrix} 0 & 0 \\ 0 & (1 - \frac{1}{m})^{|\alpha|} \end{pmatrix} + A_1^{(m)}(Z^{\alpha}) \lambda + \cdots
   \]
   Hence,
   \[
   F_m(Z_1 Z^{\alpha}) = \left(\begin{pmatrix} 0 & 0 \\ 0 & 1 - \frac{1}{m} \end{pmatrix} + \lambda \begin{pmatrix} 0 & \frac{1}{\sqrt{m}} \\ \frac{1}{\sqrt{m}} & 0 \end{pmatrix} \right) F_m(Z^{\alpha}).
   \]
   Since $A_1^{(m)}(Z^{\alpha})$ is anti-diagonal, we get that $ \langle A_1^{(m)}(Z_1 Z^{\alpha}) e_2, e_1 \rangle = \frac{1}{\sqrt{m}} (1 - \frac{1}{m})^{|\alpha|}$. Similarly, $ \langle A_1^{(m)}(Z_2 Z^{\alpha}) e_2, e_1 \rangle = -\frac{1}{\sqrt{m}} (1 - \frac{1}{m})^{|\alpha|}$. Therefore, if $p$ is a homogenous polynomial with $\deg p = n > 0$ and we write $p = \sum_{|\alpha| = n} a_{\alpha} Z^{\alpha} = Z_1 p_1 + Z_2 p_2$, then
   \[
   \langle A_1^{(m)}(p) e_2, e_1 \rangle =  \frac{1}{\sqrt{m}} (1 - \frac{1}{m})^{n-1} \sum_{|\beta| = n-1} (a_{1 \beta} - a_{2 \beta}).
   \]
   We write for $n\in \bN$, $b_n(f) = \sum_{|\beta| = n-1} (a_{1 \beta} - a_{2 \beta})$. Now recall the operators $J_n$ and $J_n'$ from Lemma \ref{lem:coeff_diff}. Write $f = \sum_{k=0}^{\infty} f_k$, where $f_k$ are the homogenous components of the Taylor-Taylor series. We note that the series converges absolutely and uniformly on $r \overline{\fD_2}$, for all $0 < r < 1$. Therefore, for every $0 < r < 1$ and every $n \in \bN$
   \[
   f\left( r J_n, r J'_n \right)^* v_0 = \sum_{k=0}^n r^k f_k(J_n,J'_n)^* v_0 =  \overline{f(0)} v_0 + \sum_{k=1}^n r^k b_k(f) v_k.
   \]
   Taking the norms, we see that $\sum_{k=1}^n r^{2k} |b_k|^2 \leq 4 \|f\|_{H^{\infty}(\fD_2)}^2$. We consider the function $B(\mu) = \sum_{k=1}^{\infty} b_k \mu^k$. Passing to a limit as $n \to \infty$ shows that $B(\mu) \in H^2(\bD)$. Now we note that for every $|\lambda| < \frac{1}{2 \sqrt{m}}$,
   \[
   \|Y_m(\lambda)\| \leq \left\| \begin{pmatrix} 0 & 0 \\ 0 & 1 - \frac{1}{m} \end{pmatrix} \right\| + \left\| \begin{pmatrix} 0 & \frac{|\lambda|}{\sqrt{m}} \\ \frac{|\lambda|}{\sqrt{m}} & 0 \end{pmatrix} \right\| < 1 - \frac{1}{2m}.
   \]
   Since the Taylor-Taylor expansion of $f$ converges absolutely and uniformly on this bidisc, we have that
   \[
   F_m(f)(\lambda) = \sum_{k=0}^{\infty} f_k(Y_m(\lambda)).
   \]
   Taking the derivative at the origin, we see that $A_1^{(m)}(f) = \sum_{k=1}^{\infty} A_1^{(m)}(f_k)$. We conclude that
   \[
   \langle A_1^{(m)}(f) e_2, e_1 \rangle = \frac{1}{(1 - \frac{1}{m})\sqrt{m}} \sum_{k=1}^{\infty} b_k(f) (1 - \frac{1}{m})^k
   \]
   For every $m \geq 2$, consider the vector $u_m \in \ell^2$ given by $u_m(k) = \frac{1}{(1 - \frac{1}{m})\sqrt{m}} (1 - \frac{1}{m})^k$. Note that
   \[
   \|u_m\|_{\ell^2}^2 = \frac{1}{m (1- \frac{1}{m})^2} \sum_{k=1}^{\infty} \left(1 - \frac{1}{m}\right)^{2k} = \frac{1}{m} \frac{1}{1 - (1 - \frac{1}{m})^2} = \frac{1}{2 - \frac{1}{m}} \leq 1.
   \]
   Now let $u_{m_k}$ be a weakly convergent subsequence with weak limit $u$. Then, for every $\ell \in \bN$
   \[
   u(\ell) = \lim_{k \to \infty} u_{m_k}(\ell) = \lim_{k\to\infty} \frac{1}{(1 - \frac{1}{m_k})\sqrt{m_k}} (1 - \frac{1}{m_k})^{\ell} = 0.
   \]
   Since every accumulation point is $0$, we conclude that $u_m$ converges to $0$ weakly. In particular,
   \[
   \lim_{m \to \infty} \langle A_1^{(m)}(f) e_2, e_1 \rangle = \lim_{m \to \infty} \langle b(f), u_m \rangle = 0.
   \]
   A similar argument (with right multiplication) shows that $\lim_{m \to \infty} \langle A_1^{(m)} e_1, e_2 \rangle = 0$. Let $d$ be a weak* accumulation point of the $d_m$. Then, Equation \eqref{eq:d_m_product} together with the preceding calculation shows that $d$ is a derivation as desired.

    It remains to define our homomorphism $\Phi \colon H^{\infty}(\fD_2) \to M_2$ given by
    \[
    \Phi(f) = \begin{pmatrix} f(0,0) & d(f) \\ 0 & f(0,0) \end{pmatrix}.
    \]
    This is clearly a continuous homomorphism. The only thing that we need to verify still is that $d \neq 0$. 
    We fix
    \begin{equation}\label{eq:f}
    f(Z_1,Z_2) = \frac{Z_1 + Z_2}{2} + \frac{Z_1 - Z_2}{2} \left(1 - \frac{Z_1 + Z_2}{2}\right)^{-1} \frac{Z_1 - Z_2}{2} .
    \end{equation}
    We need to show that $f$ is in $H^\infty(\fD_2)$; let's assume this for now and come back to it later.

    Note that
    \[
    \frac{\varphi_m(\lambda) + \varphi_m(-\lambda)}{2} = \begin{pmatrix} 0 & 0 \\ 0 & 1 - \frac{1}{m} \end{pmatrix} \text{ and } \frac{\varphi_m(\lambda) - \varphi_m(-\lambda)}{2} = \begin{pmatrix} 0 & \frac{\lambda}{\sqrt{m}} \\ \frac{\lambda}{\sqrt{m}} & 0 \end{pmatrix}.
    \]
    Hence,
    \[
    f(Y_m(\lambda)) = \begin{pmatrix} 0 & 0 \\ 0 & 1 - \frac{1}{m} \end{pmatrix} + \begin{pmatrix} 0 & \frac{\lambda}{\sqrt{m}} \\ \frac{\lambda}{\sqrt{m}} & 0 \end{pmatrix} \begin{pmatrix} 1 & 0 \\ 0 & m \end{pmatrix} \begin{pmatrix} 0 & \frac{\lambda}{\sqrt{m}} \\ \frac{\lambda}{\sqrt{m}} & 0 \end{pmatrix} = \begin{pmatrix} \lambda^2 & 0 \\ 0 & 1 - \frac{1}{m} + \frac{\lambda^2}{m} \end{pmatrix}.
    \]
    Therefore, $d_m(f) = 1$ for every $m \in \bN$. We conclude that $d(f) = 1 \neq 0$. 

    Let us now spell out the conclusion. 
    We have that $d_m(Z_j) = 0$, therefore $d(Z_j) = 0$, and thus $\Phi(Z_j) = 0$ for $j=1,2$. 
    Therefore $\Phi \in \pi^{-1}(0)$, i.e. $\Phi$ is in the fiber over the origin $0 = (0,0)$, and $\Phi$ was point evaluation, that point would have to be the origin $0 = (0,0)$. 
    But $f(0,0) = 0$, while
    \[
    \Phi(f) = \begin{pmatrix} f(0,0) & d(f) \\ 0 & f(0,0) \end{pmatrix} = \begin{pmatrix} 0 & 1 \\ 0 & 0 \end{pmatrix} .
    \]
    It therefore follows that $\Phi$ is not a point evaluation in the fiber over an interior point. 
    It also follows that $\Phi$ cannot be weak-$*$ continuous. 
\end{example} 

\noindent{\em Showing that $f$ is a bounded NC function.} 
Recall that $f$ is given by \eqref{eq:f}. 
The fact that $f \in H^\infty(\fD_2)$, in fact $\|f\|\leq 1$, follows from standard realization techniques, see \cite{agler2015} or in particular Definition 1.1 and Theorem 1.2 in \cite{agler2018global}. 
Since we need only the easier part of that result, we can give a brief and self-contained proof. 

Fix $(X_1, X_2) \in \fD_2$, say $X_1$ and $X_2$ are contractions in $M_n$. 
We will show the inequality $\|f(X_1,X_2)u\| \leq \|u\|$ for all $u \in \bC^n$. 
Put 
\[
Y_1 = \frac{X_1+X_2}{2} \,\, \textrm{ and } \,\, Y_2 = \frac{X_1 - X_2}{2}.
\]
The matrix 
\[
T = \begin{pmatrix}
Y_1 & Y_2 \\ Y_2 & Y_1
\end{pmatrix}
\]
is unitarily equivalent to $\diag(X_1,X_2)$, and so $\|T\|<1$. 
Now, given $u \in \bC^n$, define $v = (I-Y_1)^{-1}Y_2 u$. 
Then a direct calculation shows that 
\[
T \begin{pmatrix} u \\ v \end{pmatrix} =
\begin{pmatrix} f(X_1,X_2) u \\ v \end{pmatrix}, 
\]
and since $\|T\|<1$ we get 
\[
\|f(X_1,X_2)u\|^2 + \|v\|^2 \leq \|u\|^2 + \|v\|^2, 
\]
therefore $\|f(X_1,X_2)u\| \leq \|u\|$, as claimed.

\subsection{The free pseudohyperbolic metrics}

We now define metrics on the representation spaces $\operatorname{Rep}_{n}(\cA(\fV))$ that will give the varieties the added structure they need to become complete invariants for the algebras. 
For $\Phi, \Psi \in \operatorname{Rep}_{n}(\cA(\fV))$, define
\[
\delta_{cb}(\Phi, \Psi) = \|\Phi - \Psi\|_{cb}
\]
and for $X, Y \in \widetilde{\ol{\fV}}(n)$
\begin{equation}\label{eq:delta_cb}
\delta_{cb}(X, Y) = \|\Phi_X - \Phi_Y\|_{cb}. 
\end{equation}
Similarly, we define
\[
\delta_b(\Phi, \Psi) = \|\Phi - \Psi\|
\]
and
\begin{equation}\label{eq:delta_b}
\delta_{b}(X, Y) = \|\Phi_X - \Phi_Y\|. 
\end{equation}

\begin{proposition}\label{eq:psd_contractive}
If $G \colon \ol{\fW} \to \ol{\fV}$ is an NC map with coordinates in $\cA(\fW)$ and $\kappa \in \{b,cb\}$, then for all $X,Y \in \ol{\fW}$
\[
\delta_\kappa (G(X),G(Y)) \leq \delta_\kappa (X,Y) . 
\]
The same conclusion holds for NC maps between the non-closed varieties $\fW$ and $\fV$. 
\end{proposition}
\begin{proof}
If $F$ is in the closed unit ball $\cA(\fV)_1$ of $\cA(\fV)$, then $F \circ G$ is in the closed unit ball of $\cA(\fW)$, thus
\[
\sup_{W\in \fW} \|F(G(W))\| \leq \sup_{V\in \fV} \|F(V)\| .
\]
By the definition of $\delta_b$, 
\begin{align*}
\delta_b(G(X),G(Y)) 
&= \|\Phi_{G(X)} - \Phi_{G(Y)}\| \\
&= \sup_{F \in \cA(\fV)_1} \|F(G(X)) - F(G(Y))\| \\
&\leq \sup_{H \in \cA(\fW)_1} \|H(X) - H(Y)\| \\
&= \delta_b(X,Y). 
\end{align*}
The same argument works for $\kappa = cb$. 
\end{proof}

Let $\fV \subseteq \bB_E$ be a homogeneous NC variety.  
For $X \in \widetilde{\fV}$, the evaluation representation $\Phi_X$ can be considered as a representation of $\cA(\fV) = A(\fV)$  or as a representation of $\cA(\fV) = H^\infty(\fV)$.
Thus, one can define distance functions on $\widetilde{\fV}$ using equations \eqref{eq:delta_cb} and \eqref{eq:delta_b}, thinking of $\Phi_X,\Phi_Y$ as elements of $\operatorname{Rep}_{n}(A(\fV))$ or $\operatorname{Rep}_{n}(H^\infty(\fV))$; one can also define distances thinking of $\Phi_X,\Phi_Y$ as elements of $\operatorname{Rep}_{n}(A(\bB_E))$ or $\operatorname{Rep}_{n}(H^\infty(\bB_E))$.
It turns out that all the reasonable notions agree. 

\begin{proposition}\label{prop:equi_F_norm}
Let $\fV\subseteq\bB_E$ be an NC variety. Then, for every
$n\geq1$ and every $X,Y\in\widetilde{\fV}(n)$,
\begin{align*}
\sup_{F\in H^\infty(\fV)_1}
\|\widetilde F(X)-\widetilde F(Y)\|
&=
\sup_{F\in H^\infty(\bB_E)_1}
\|\widetilde F(X)-\widetilde F(Y)\| \\
&=
\sup_{F\in A(\bB_E)_1}
\|\widetilde F(X)-\widetilde F(Y)\| \\
&= 
\sup_{F\in A(\fV)_1}
\|\widetilde F(X)-\widetilde F(Y)\|. 
\end{align*}
\end{proposition}
\begin{proof}
Let $F\in H^\infty(\bB_E)_1$. Then
$F|_{\fV}\in H^\infty(\fV)_1$. 
Moreover, for
$X=S^{-1}X_0S\in\widetilde{\fV}$, where
$X_0\in\fV$ and $S$ is an invertible matrix,
$\widetilde F(X)
=
\widetilde{F|_{\fV}}(X).$
Thus
\[
\sup_{F\in H^\infty(\bB_E)_1}
\|\widetilde F(X)-\widetilde F(Y)\|
\leq
\sup_{F\in H^\infty(\fV)_1}
\|\widetilde F(X)-\widetilde F(Y)\|.
\]
For the reverse inequality, let
$F\in H^\infty(\fV)_1$. 
By \cite[Corollary 3.4]{BMV-int}
there exists
$\widehat F\in H^\infty(\bB_E)$
such that
$\widehat F|_{\fV}=F$
and
$\|\widehat F\| = \|F\|$.
Moreover, for
$X\in\widetilde{\fV}$ 
we have $\widetilde{\widehat F}(X)=\widetilde F(X)$, 
hence
\[
\sup_{F\in H^\infty(\fV)_1}
\|\widetilde F(X)-\widetilde F(Y)\|
=
\sup_{F\in H^\infty(\bB_E)_1}
\|\widetilde F(X)-\widetilde F(Y)\|.
\]
We next show that
\[
\sup_{F\in H^\infty(\bB_E)_1}
\|\widetilde F(X)-\widetilde F(Y)\|
=
\sup_{F\in A(\bB_E)_1}
\|\widetilde F(X)-\widetilde F(Y)\|.
\]
Since
$A(\bB_E)_1\subseteq H^\infty(\bB_E)_1$,
we have
\[
\sup_{F\in A(\bB_E)_1}
\|\widetilde F(X)-\widetilde F(Y)\|
\leq
\sup_{F\in H^\infty(\bB_E)_1}
\|\widetilde F(X)-\widetilde F(Y)\|.
\]
For the converse, let
$F\in H^\infty(\bB_E)_1$
and, for $0<r<1$, define
$F_r(Z)=F(rZ),
\, Z\in\bB_E.$
Since $\bB_E$ is homogeneous,
\[
\|F_r\|
=
\sup_{Z\in\bB_E}\|F(rZ)\|
\leq
\|F\|
\leq1.
\]
By radial approximation,
$F_r\in A(\bB_E)_1$
and 
$\widetilde{F_r}(X)\longrightarrow\widetilde F(X)$
as $r\to1^-$. 
Consequently,
\[
\|
\widetilde{F_r}(X)-\widetilde{F_r}(Y)
\|
\longrightarrow
\|
\widetilde F(X)-\widetilde F(Y)
\|.
\]
Since $F_r\in A(\bB_E)_1$, we obtain
\[
\|\widetilde F(X)-\widetilde F(Y)\|
\leq
\sup_{H\in A(\bB_E)_1}
\|\widetilde H(X)-\widetilde H(Y)\|, 
\]
and thus
\[
\sup_{F\in H^\infty(\bB_E)_1}
\|\widetilde F(X)-\widetilde F(Y)\|
=
\sup_{F\in A(\bB_E)_1}
\|\widetilde F(X)-\widetilde F(Y)\|.
\]

Finally, the restriction map
$
R_{\fV}\colon A(\bB_E)\longrightarrow A(\fV)
$
is well defined and contractive. Indeed, if
$F\in A(\bB_E)$ and $(p_j)$ is a sequence of NC polynomials
converging to $F$ in $A(\bB_E)$, then
\[
\|p_j|_{\fV}-F|_{\fV}\|_{A(\fV)}
\leq
\|p_j-F\|_{A(\bB_E)}
\longrightarrow0,
\]
so $F|_{\fV}\in A(\fV)$. Moreover, since
$\fV\subseteq\bB_E$,
$\|F|_{\fV}\|_{A(\fV)}
=
\sup_{Z\in\fV}\|F(Z)\|
\leq
\sup_{Z\in\bB_E}\|F(Z)\|
=
\|F\|_{A(\bB_E)}.$
Also,
$
\widetilde{F|_{\fV}}(Z)=\widetilde F(Z)
$
for
$
Z\in\widetilde{\fV}.
$
Consequently,
\[
\sup_{F\in A(\bB_E)_1}
\|\widetilde F(X)-\widetilde F(Y)\|
\leq
\sup_{F\in A(\fV)_1}
\|\widetilde F(X)-\widetilde F(Y)\|.
\]
On the other hand, since
$A(\fV)_1\subseteq H^\infty(\fV)_1,$
we have
\[
\sup_{F\in A(\fV)_1}
\|\widetilde F(X)-\widetilde F(Y)\|
\leq
\sup_{F\in H^\infty(\fV)_1}
\|\widetilde F(X)-\widetilde F(Y)\|.
\]
Combining the preceding equalities and inequalities, we obtain
\[
\begin{aligned}
\sup_{F\in H^\infty(\fV)_1}
\|\widetilde F(X)-\widetilde F(Y)\|
&=
\sup_{F\in H^\infty(\bB_E)_1}
\|\widetilde F(X)-\widetilde F(Y)\|=
\sup_{F\in A(\bB_E)_1}
\|\widetilde F(X)-\widetilde F(Y)\|\\
&\leq
\sup_{F\in A(\fV)_1}
\|\widetilde F(X)-\widetilde F(Y)\| \leq
\sup_{F\in H^\infty(\fV)_1}
\|\widetilde F(X)-\widetilde F(Y)\|.
\end{aligned}
\]
Therefore, all the inequalities are equalities, which proves the result.
\end{proof}

\begin{remark}
We see that the metric $\delta_{b}$ is independent of whether it is defined using
$H^\infty(\bB_E)$ or $A(\bB_E)$ or $H^\infty(\fV)$ (or $A(\fV)$ when applicable)  in the definition \eqref{eq:delta_b}.    
\end{remark}

\begin{remark}
For every $k\geq1$, the same argument at the matrix level shows that
the pseudohyperbolic
metric $\delta_{cb}$ is independent of which algebra is taken in \eqref{eq:delta_cb}.
\end{remark}

\section{Classification of algebras of uniformly continuous NC functions}

Throughout this section, we let $E$ and $F$ be two finite-dimensional operator space structures and let $\fV \subseteq \bB_E$ and $\fW \subseteq \bB_F$ be two homogeneous NC varieties. 
Our goal is to obtain necessary and sufficient conditions for the existence of a (completely) bounded isomorphism $\alpha \colon A(\fV) \to A(\fW)$ in terms of the geometry of the similarity envelopes of the varieties, that is, in terms of the existence of certain maps $G \colon \widetilde{\fW} \to \widetilde{\fV}$. 
Recall that if $\alpha$ is a bounded isomorphism between Banach algebras, then $\alpha^{-1}$ is bounded, by the inverse mapping theorem. 
In the completely bounded category there is no such results, thus, throughout, when we say that a map $\alpha \colon H^\infty(\fV) \to H^\infty(\fW)$ is a {\em completely bounded isomorphism}, we mean that $\alpha$ is an isomorphism of algebras such that both $\alpha$ and its inverse $\alpha^{-1}$ are completely bounded.

The following lemma will be useful. 

\begin{lemma}\label{lem:uniform_components}
Suppose that $\fV\subseteq\bB_E$ is a variety and that
$\fW\subseteq\bB_F$ is a homogeneous variety. Let
$G\colon\widetilde{\fW}\longrightarrow\widetilde{\fV}$
be an NC holomorphic map such that
$G=(G_1,\ldots,G_d)$. Then the following are equivalent:
\begin{enumerate}
    \item The restriction
   $G|_{\fW}\colon\fW\longrightarrow\widetilde{\fV}$
    is uniformly continuous with respect to the ambient operator-space
    norms.

    \item
    $G_j|_{\fW}\in A(\fW),
    $ for $ j=1,\ldots,d.$
\end{enumerate}
\end{lemma}
\subsection{Maps induce homomorphisms}

The following results have straightforward proofs, similar to the those in the case of the row ball (see \cite[Propositions 5.2 and 5.5]{SSS2}).

\begin{proposition}\label{prop:equiv_G_norm}
Suppose that $\fV \subseteq \bB_E$ is a variety and $\fW \subseteq \bB_F$ is a homogeneous varieties, and let $G \colon \widetilde{{\fW}} \to \widetilde{{\fV}}$ an NC holomorphic function such that
$G|_{\fW}:\fW\longrightarrow\widetilde{\fV}$
is uniformly continuous. Let $\kappa \in \{b, cb\}$. Then the following statements are equivalent:
\begin{enumerate}
\item $G$ is Lipschitz with respect to $\delta_{\kappa}$,
\item $\sup_{W\in \fW}\|\Phi_{G(W)}\|_{\kappa} < \infty$,
\item the mapping $f \mapsto \widetilde{f} \circ G|_{\fW}$ for $f \in A(\fV)$ is a well-defined unital bounded homomorphism $\alpha:A(\fV) \to A(\fW)$ which is completely bounded in the case $\kappa = cb$.
\end{enumerate}
Furthermore, in this case
\begin{itemize}
\item[(a)] $\|\alpha \|_{\kappa} = \sup_{W\in \fW}\|\Phi_{G(W)}\|_{\kappa}$,
\item[(b)] the value in (a) is a Lipschitz constant for $G$. 
\end{itemize}
\end{proposition}

\begin{proof}
We prove
$(1)\Longrightarrow(2)\Longrightarrow(3)\Longrightarrow(1).$

\medskip
\noindent
\emph{$(1)\Rightarrow(2)$.}
Suppose that $G$ is $\delta_\kappa$-Lipschitz with constant $L$.
Since $\fW$ is homogeneous,
$0_n\in\fW(n),$ for $ n\geq1.$
Let $W\in\fW(n)$. Then
\begin{align*}
\|\Phi_{G(W)}\|_\kappa
&\leq
\|\Phi_{G(W)}-\Phi_{G(0_n)}\|_\kappa
+
\|\Phi_{G(0_n)}\|_\kappa \\
&=
\delta_\kappa\bigl(G(W),G(0_n)\bigr)
+
\|\Phi_{G(0_n)}\|_\kappa \\
& \leq
L\,\delta_\kappa(W,0_n)
+
\|\Phi_{G(0_n)}\|_\kappa.
\end{align*}
Since $W,0_n\in\fW(n)$, the evaluation representations
$\Phi_W$ and $\Phi_{0_n}$ are completely contractive. Hence
\[
\delta_\kappa(W,0_n)
=
\|\Phi_W-\Phi_{0_n}\|_\kappa
\leq
\|\Phi_W\|_\kappa+\|\Phi_{0_n}\|_\kappa
\leq2.
\]
Moreover, since $G$ is an NC map,
\[
G(0_n)
=
\underbrace{
G(0_1)\oplus\cdots\oplus G(0_1)
}_{n\text{ times}}.
\]
Consequently,
$\|\Phi_{G(0_n)}\|_\kappa
=
\|\Phi_{G(0_1)}\|_\kappa
$ for $
\kappa\in\{b,cb\}.$
Since
$G(0_1)\in \widetilde{\fV}(1)=\fV(1),$
the evaluation representation
$\Phi_{G(0_1)}$
is completely contractive.
Since $\Phi_{G(0_1)}$ is unital, 
$\|\Phi_{G(0_1)}\|_{\kappa}
=
1.$
Therefore,
\[
\|\Phi_{G(W)}\|_\kappa
\leq 2L+1,
\qquad W\in\fW,
\]
and hence
$\sup_{W\in\fW}\|\Phi_{G(W)}\|_\kappa<\infty.$

\medskip
\noindent
\emph{$(2)\Rightarrow(3)$.}
Set
$M_G
:=
\sup_{W\in\fW}\|\Phi_{G(W)}\|_\kappa
<\infty.$
We show that
$\widetilde f\circ G \in A(\fW)$ for every $f\in A(\fV)$
 
Fix $f\in A(\fV)$. There exists a sequence of NC
polynomials $(p_k)_{k=1}^\infty$ such that
$\|p_k|_{\fV}-f\|_{A(\fV)}
\longrightarrow 0$.
By our blanket assumption on $G$ and Lemma \ref{lem:uniform_components},
$G_j|_{\fW}\in A(\fW),$ for $ j=1,\ldots,d$, therefore
$p_k \circ G|_{\fW}
=
p_k(G_1|_{\fW},\ldots,G_d|_{\fW})
\in A(\fW)$ for all $k$. 
For $k, m\geq1$,
\begin{align*}
\|p_k\circ G|_{\fW}-p_m\circ G|_{\fW}\|_{A(\fW)}
&=
\sup_{W\in\fW}
\|\Phi_{G(W)}(p_k-p_m)\| \\
&\leq
M_G\,
\|p_k-p_m\|_{A(\fV)}.
\end{align*}
Thus $(p_k\circ G|_{\fW})$ is Cauchy in $A(\fW)$, so there exists $h\in A(\fW)$ such that
$p_k\circ G|_{\fW}
\longrightarrow h \in A(\fW)$. 
It is now not hard to see that
\[
\widetilde f\circ G = h
\in A(\fW).
\]
Consequently, the map
$\alpha:A(\fV)\longrightarrow A(\fW),
$ given by $
\alpha(f)
=
\left.\widetilde f\circ G\right|_{\fW}$
is well defined.
It is immediate that $\alpha$ is a unital algebra homomorphism.

Now suppose that $\kappa=cb$ ($\kappa = b$ is similar). Let
$[f_{ij}]\in M_k(A(\fV)).$
Then
\begin{align*}
\|\alpha^{(k)}([f_{ij}])\|_{M_k(A(\fW))}
&=
\sup_{W\in\fW}
\left\|
[\widetilde f_{ij}(G(W))]
\right\|=
\sup_{W\in\fW}
\left\|
\Phi_{G(W)}^{(k)}([f_{ij}])
\right\|\leq
M_G
\|[f_{ij}]\|_{M_k(A(\fV))}.
\end{align*}
Taking the supremum over $k\geq1$, we obtain
$
\|\alpha\|_{cb}\leq M_G.$

Conversely, for every $W\in\fW$ and $f\in A(\fV)$,
\[
\Phi_{G(W)}(f)
=
\widetilde f(G(W))
=
\alpha(f)(W)
=
\Phi_W(\alpha(f)).
\]
Hence
$\Phi_{G(W)}
=
\Phi_W \circ \alpha$. 
Since $\Phi_W$ is completely contractive,
$\|\Phi_{G(W)}\|_\kappa
\leq
\|\alpha\|_\kappa.$
Taking the supremum over $W\in\fW$, we obtain
$M_G
\leq
\|\alpha\|_\kappa.$
Combining the two inequalities gives {\rm(a)}.

\medskip
\noindent
\emph{$(3)\Rightarrow(1)$.}
Let
$X,Y\in\widetilde{\fW}(n).$
We first show that
$\Phi_{G(X)}
=
\Phi_X\circ\alpha.$

Since $X\in\widetilde{\fW}(n)$, there exist $W\in\fW(n)$ and
$S\in\GL_n(\mathbb C)$ such that
$X=SWS^{-1}.$
Hence, for every $f\in A(\fV)$,
\begin{align*}
\widetilde f(G(X))
=
\widetilde f\bigl(SG(W)S^{-1}\bigr)=
S\widetilde f(G(W))S^{-1}=
S\alpha(f)(W)S^{-1}=
\widetilde{\alpha(f)}(X).
\end{align*}
Therefore
$\Phi_{G(X)}(f)
=
(\Phi_X\circ\alpha)(f),
$ for $
f\in A(\fV),$
and hence
$\Phi_{G(X)}
=
\Phi_X\circ\alpha.$
Similarly,
$\Phi_{G(Y)}
=
\Phi_Y\circ\alpha.$
Consequently,
\begin{align*}
\delta_\kappa(G(X),G(Y))
&=
\|\Phi_{G(X)}-\Phi_{G(Y)}\|_\kappa \\
& =
\|(\Phi_X-\Phi_Y)\circ\alpha\|_\kappa \\
& \leq
\|\alpha\|_\kappa\,
\|\Phi_X-\Phi_Y\|_\kappa =
\|\alpha\|_\kappa\,
\delta_\kappa(X,Y).
\end{align*}
Thus $G$ is $\delta_\kappa$-Lipschitz.
The value in {\rm(a)} is a Lipschitz constant for $G$. This proves
{\rm(b)} and completes the proof.
\end{proof}

\begin{remark}\label{rem:replace_closure}
In the above proposition the varieties $\widetilde{\fV}, \widetilde{\fW}$ can be replaced by the similarity envelopes of the closures $\widetilde{\ol{\fV}}, \widetilde{\ol{\fW}}$. 
\end{remark}

\subsection{Homomorphisms induce maps}

Any homomorphism $\alpha \colon \cA \to \cB$ between algebras gives rise to a map between representation spaces
\[
\alpha^* \colon \operatorname{Rep}_{\operatorname{fin}}(\cB) \to \operatorname{Rep}_{\operatorname{fin}}(\cA)
\]
given by 
\[
\alpha^* \Phi = \Phi \circ \alpha. 
\]
If $\alpha$ is (completely) bounded then $\alpha^*$ maps (completely) bounded representations to (completely) bounded representation. 
It is plain to see that $\alpha^*$ preserves grading and intertwines similarities. 

By the identification of the representation spaces of algebras of uniformly continuous functions given in Proposition \ref{prop:rep_AV}, we shall obtain in Proposition \ref{prop:homo_ind_NCmap} below a necessary condition for the existence of a bounded homomorphism $\alpha \colon A(\fV) \to A(\fW)$. 

\begin{proposition}\label{prop:homo_ind_NCmap}
    Suppose that $\fV \subseteq \bB_E$ and $\fW \subseteq \bB_F$ are homogeneous varieties, and let $\alpha \colon A(\fV) \to A(\fW)$ be a bounded isomorphism. Then there exists an NC bijection $G \colon \widetilde{\ol{\fW}} \to \widetilde{\ol{\fV}}$ whose restriction to $\fW$ has components in $A(\fW)$, that restricts to an NC bijection $G \colon \widetilde{\fW} \to \widetilde{\fV}$ that is $\delta_b$-bi-Lipschitz ($\delta_{cb}$-bi-Lipschitz if $\alpha$ is a completely bounded isomorphism), 
    and $\alpha$ is given by 
    \[
    \alpha(f) = \widetilde{f}\circ G|_{{\fW}}, \qquad \forall \, f\in A(\fV).
    \]
\end{proposition}

\begin{proof}
Given $W\in\widetilde{\ol{\fW}}$, Proposition \ref{prop:rep_AV} gives us $\Phi_W\in\operatorname{Rep}_{\mathrm{fin}}(A(\fW))$.
Since $\alpha:A(\fV)\to A(\fW)$ is bounded,
$
\Phi_W\circ\alpha \in\operatorname{Rep}_{\mathrm{fin}}(A(\fV))$.
Hence, again by Proposition \ref{prop:rep_AV}, there exists
a unique
$G(W)\in\widetilde{\ol{\fV}}$
such that
\begin{equation}\label{eq:def-induced-G}
\Phi_{G(W)}
=
\Phi_W\circ\alpha.
\end{equation}
Thus we obtain a well-defined map
$G:\widetilde{\ol{\fW}}
\longrightarrow
\widetilde{\ol{\fV}}.$
Let
$Z_1,\ldots,Z_d$ denote the coordinate functions on $\fV$.
For $j=1,\ldots,d$ and
$W\in\widetilde{\ol{\fW}}$, by
\eqref{eq:def-induced-G},
\[
\begin{aligned}
G_j(W)
=
\Phi_{G(W)}(Z_j)=
(\Phi_W\circ\alpha)(Z_j)=
\widetilde{\alpha(Z_j)}(W), 
\end{aligned}
\]
hence $G$ is an NC map and $G\big|_\fW$ has components in $A(\fW)$, as claimed.

Apply the same argument to $\alpha^{-1}:A(\fW)\longrightarrow A(\fV)$. 
Let $ H \colon \widetilde{\ol{\fV}} \longrightarrow \widetilde{\ol{\fW}}$ be the NC map induced by $\alpha^{-1}$. 
For $X\in\widetilde{\ol{\fW}}$, 
\[
\Phi_{H(G(X))} = \Phi_{G(X)}\circ\alpha^{-1} = (\Phi_X\circ\alpha)\circ\alpha^{-1} 
= \Phi_X.
\]
By injectivity of the correspondence $X\leftrightarrow\Phi_X$, $H(G(X))=X$. 
Similarly, $G(H(Y))=Y$ for all $Y\in\widetilde{\ol{\fV}}$. 
Thus $H=G^{-1}$.

We next prove that $G(\widetilde{\fW}) = \widetilde{\fV}$. 
By symmetry, it suffices to prove that $G(\widetilde{\fW}) \subseteq \widetilde{\fV}$. 
By similarity invariance, it suffices to show that $G(\fW) \subseteq \widetilde{\fV}$. 
Fix $X\in\fW(n)$. 
Suppose, toward a contradiction, that $G(X)\notin\widetilde{\fV}(n)$. 
Then $\rho_E(G(X))=1$. 
Choose $r>1$ such that $zrX\in\fW(n)$ for $z \in \mathbb D.$ Define the holomorphic map
    \[
    \Psi:\mathbb D
    \longrightarrow
    \widetilde{\ol{\fV}}(n),
    \qquad
    \Psi(z):=G(zrX).
    \]
    If
    $\Psi(z)\in\widetilde{\fV}(n)$,
    then
    $\rho_E(\Psi(z))<1,$
    since $\Psi(z)$ is similar to a point of $\fV(n)\subseteq\bB_E(n)$.
    Likewise, if
    $\Psi(z)\in
    \widetilde{\ol{\fV}}(n)
    \setminus
    \widetilde{\fV}(n),$
    then $\rho_E(\Psi(z))=1$. 
    Thus
    \[
    \rho_E(\Psi(z))\leq1,
    \qquad
    z\in\mathbb D.
    \] 
    Since $1/r\in\mathbb D$,
    \[
    \rho_E\left(
    \Psi\left(\frac1r\right)
    \right)
    =
    \rho_E(G(X))
    =
    1.
    \]
    By
    \cite[Lemma 2.5]{shalsham2025spectral},
    the function
    $ z\longmapsto\rho_E(\Psi(z))$
    is log-subharmonic. Since it is bounded above by $1$ and
    attains the value $1$ at the interior point $1/r$, the maximum
    principle gives
    \[
    \rho_E(\Psi(z))=1,
    \qquad
    z\in\mathbb D.
    \]
    In particular, $\rho_E(G(0_n))=1$. Applying the same reasoning to any complex disc
in $\fW(n)$ passing through the origin, we find that
$\rho_E(G(Y))=1$ for every $Y\in\fW(n)$.
Therefore,
$G(\fW(n))
\subseteq
\widetilde{\ol{\fV}}(n)\setminus\widetilde{\fV}(n).$
By continuity, upper semicontinuity of $\rho_E$,
and similarity invariance, we obtain
\[
G(\widetilde{\ol{\fW}}(n))
\subseteq
\widetilde{\ol{\fV}}(n)\setminus\widetilde{\fV}(n).
\]
Since $G$ is surjective onto
$\widetilde{\ol{\fV}}$, this is a contradiction,
because $0_n\in\widetilde{\fV}(n)$.
Therefore $G(\fW)\subseteq\widetilde{\fV}$,
and hence
$G(\widetilde{\fW})=\widetilde{\fV}$.

Finally, we prove that $G$ bi-Lipschitz. 
We have $\Phi_{G(W)} = \Phi_W\circ\alpha$ for $W\in\fW$. 
Since $\Phi_W$ is completely contractive,
 $ \sup_{W\in\fW}
    \|\Phi_{G(W)}\|
    \leq
    \|\alpha\|$, 
    hence by Proposition \ref{prop:equiv_G_norm},
    $G$ is $\delta_b$-Lipschitz. So for every $n\in\mathbb N$ and
$X,Y\in\widetilde{\fW}(n)$,
\[
\delta_b^{\fV}(G(X),G(Y))
\leq
\|\alpha\|\,
\delta_b^{\fW}(X,Y).
\]
 Applying the same argument to $\alpha^{-1}$ and $H=G^{-1}$ for every $n\in\mathbb N$ and
$U,V\in\widetilde{\fV}(n)$,
    \[
    \delta_b^{\fW}(H(U),H(V))
    \leq
    \|\alpha^{-1}\|\,
    \delta_b^{\fV}(U,V).
    \]
    Taking
    $U=G(X)$ and $V=G(Y)$ yields
    \[
    \frac{1}{\|\alpha^{-1}\|}
    \delta_b^{\fW}(X,Y)
    \leq
    \delta_b^{\fV}(G(X),G(Y))
    \leq
    \|\alpha\|
    \delta_b^{\fW}(X,Y).
    \]
    Therefore $G$ is $\delta_b$-bi-Lipschitz.

    If $\alpha$ is completely bounded, then
  $ \sup_{W\in\fW}
    \|\Phi_{G(W)}\|_{cb}
    \leq
    \|\alpha\|_{cb},$
    and Proposition \ref{prop:equiv_G_norm} implies $G$ is $\delta_{cb}$-Lipschitz.
   Applying the same argument to the completely bounded inverse $\alpha^{-1}$, we obtain, for every 
   $n\in\mathbb N$ and every $X,Y\in\widetilde{\fW}(n)$,
\[
\frac{1}{\|\alpha^{-1}\|_{cb}}
\delta_{cb}^{\fW}(X,Y)
\leq
\delta_{cb}^{\fV}(G(X),G(Y))
\leq
\|\alpha\|_{cb}\,
\delta_{cb}^{\fW}(X,Y).
\]
Hence $G$ is $\delta_{cb}$-bi-Lipschitz.

    Finally, for $f\in A(\fV)$ and $W\in\fW$,
    \[
    \alpha(f)(W)
    =
    \Phi_W(\alpha(f))
    =
    (\Phi_W\circ\alpha)(f)
    =
    \Phi_{G(W)}(f)
    =
    \widetilde f(G(W)).
    \]
    Therefore
    $ \alpha(f)
    =
    \left.\widetilde f\circ G\right|_{\fW},$
    as claimed.
\end{proof}

\subsection{Classification}

We now arrive at the basic if-and-only-if classification up to (completely) bounded isomorphisms theorem. 
In the next section we shall obtain a significant refinement. 

\begin{theorem}
\label{thm:bounded-cb-classification}
Let $\fV \subseteq \bB_E$ and $\fW \subseteq \bB_F$ be two homogeneous varieties. 
For $\kappa \in \{b,cb\}$, the following are equivalent: 
\begin{enumerate}
    \item There is an algebra isomorphism $\alpha \colon A(\fV) \to A(\fW)$ with $\|\alpha\|_\kappa< \infty$ and $\|\alpha^{-1}\|_\kappa < \infty$. 
    \item There exists a bijective and $\delta_\kappa$-bi-Lipschitz NC map $G: \widetilde{\fW} \to \widetilde{\fV}$ such that $G\big|_\fW$ and $G^{-1}\big|_\fV$ are uniformly continuous. 
\end{enumerate}
Every  isomorphism $\alpha$ as in (1) is implemented by pre-composition with a map $G$ as in (2):
\[
\alpha(f) = \widetilde{f} \circ G \qquad , \qquad f \in A(\fV).
\]
\end{theorem}

\begin{proof}
This follows from combining Propositions \ref{prop:equiv_G_norm} and \ref{prop:homo_ind_NCmap}. 
\end{proof}

In Section \ref{subsec:cnter_examples1} examples will be provided showing that a uniformly continuous map on $\fW$ need not be Lipschitz, and that a Lipschitz map need not have uniformly continuous components. 
However, we will improve the above result in Theorem \ref{thm:linearisation_maps} by showing that the {\em existence} of a bijective NC $G \colon \widetilde{\fW} \to \widetilde{\fV}$ that is $\delta_\kappa$-bi-Lipschitz implies the {\em existence} of a $\kappa$ bounded isomorphism $\alpha \colon A(\fV) \to A(\fW)$, albeit $\alpha$ need not be implemented via pre-composition with $G$. 

\section{Linearization}\label{sec:linear}

In previous works on the isomorphism problem for algebras of multipliers or algebras of bounded NC functions \cite{AL12,DRS1,Hartz2,SSS1,SSS2,SS09}, it was seen that the existence of an isomorphism between algebras on homogeneous varieties implies not only the existence of an appropriate biholomorphism between the varieties, but also the existence of a linear map between the varieties, which induces an isomorphism. 
Thus an isomorphism between such algebras implies an equivalence of a very special and rigid kind between the underlying geometrical objects. 

A key step in showing that the existence of an isomorphism implies the existence of a linear map between the varieties, is showing that the existence of an isomorphism implies the existence of a vacuum preserving isomorphism, i.e., showing the existence of an isomorphism $\beta$ whose induced map $\beta^*$ sends $\Phi_0$ to $\Phi_0$. 
In situations where the ambient domain was a Euclidean ball, or an NC row ball, one could use the so-called ``disc trick" \cite{shalit2019disc}, which shows very broadly that whenever there is a biholomorphism between two homogeneous varieties, there is also one that maps the origin to the origin. A key observation in that setting was that a biholomorphism between Euclidean balls maps a linear disc onto a linear disc. 
The disc trick then amounts to taking an arbitrary image of the origin under a biholomorphism, and rotating it using the two gauge transformations on both sides until it hits the origin. 

When the ambient domains are not Euclidean balls, a biholomorphism between balls need not map a disc onto a disc. 
One then needs new methods of finding a biholomorphism mapping the origin to the origin. 
In \cite{SampatShalit25} an NC version of a theorem of Braun, Kaup and Upmeier \cite{BKU78} was used as a replacement. 
However, this method does not generalize to the current setting, where the NC varieties and domains of interest are similarity envelopes, hence are unbounded. 

The following Lemma is the key to solving the problem in the generality studied here. 

\begin{lemma}[Two-circle trick]\label{lem:two_circ}
Let $D\subseteq \mathbb{C}^{d}$ be a bounded convex domain, and let $K$ and $L$ be connected compact subgroups of $\operatorname{Aut}(D)$. Define
\[
\operatorname{Fix}(K)
=
\{z\in D : k(z)=z \text{ for every } k\in K\},
\]
and similarly for $L$. Suppose that
\[
\operatorname{Fix}(K)=\{0\},
\qquad
\operatorname{Fix}(L)=\{p\}.
\]
Then there are $k_j\in K$ and $\ell_j\in L$, $1\leq j\leq m$,
such that
\[
w=k_1\ell_1\cdots k_m\ell_m
\]
satisfies
\[
w(p)=0.
\]
\end{lemma}

\begin{proof}
Let
\[
G=\operatorname{Aut}(D)
\]
with the compact-open topology. By Cartan's theorem, $G$ has a natural
structure of a finite-dimensional real Lie group, compatible with the
compact-open topology, and its action on $D$ is real analytic; see
\cite[Chapter~IX, p.~169]{Narasimhan1971}.

Since $G$ is Hausdorff and $K$ and $L$ are compact, they are closed
subgroups of $G$. By the closed subgroup theorem, $K$ and $L$ are
embedded Lie subgroups of $G$; see
\cite[Theorem~9.3.7]{HilgertNeeb2012}. Since they are connected
manifolds, they are path connected.

Let
\[
H=\langle K,L\rangle
\]
be the abstract subgroup of $G$ generated by $K$ and $L$. No
topological closure is taken here. Since the identity belongs to both
$K$ and $L$, we have
\begin{equation}\label{eq:Hunion}
H=\bigcup_{m\geq 1}(KL)^m.
\end{equation}
Indeed, consecutive factors belonging to the same subgroup may be
multiplied together, and identity factors may be inserted to make a word
alternate between $K$ and $L$.

For every $m$, the set $(KL)^m$ is path connected, since it is the
continuous image of the path-connected space $(K\times L)^m$ under
the multiplication map. The sets in \eqref{eq:Hunion} form an increasing family and
all contain the identity. Hence $H$ is path connected. In particular,
\[
H\subseteq G^{\circ},
\]
where $G^{\circ}$ denotes the identity component of $G$.

By Yamabe's theorem, every path-connected subgroup of a connected Lie
group is an integral Lie subgroup; see
\cite[Theorem~9.6.1]{HilgertNeeb2012}. Moreover, the initial subgroup
theorem equips $H$ with its canonical connected Lie group structure
such that the inclusion
\[
    \iota:H\longrightarrow G^{\circ}
\]
is an initial Lie subgroup; see
\cite[Definition~9.6.11 and Theorem~9.6.13]{HilgertNeeb2012}.

We record a topological point that will be needed below. The inclusion
$K\to G^{\circ}$ is smooth and has image contained in $H$. By the
initial property of $H$, its corestriction $K\to H$ is smooth and
therefore continuous. Thus $K$ is compact when regarded as a subgroup
of the intrinsic Lie group $H$. The same argument applies to $L$.

The maximal compact subgroup theorem states that, in a Lie group with
finitely many connected components, every compact subgroup is contained
in a maximal compact subgroup and all maximal compact subgroups are
conjugate; see \cite[Theorem~14.1.3]{HilgertNeeb2012}. Since $H$ is
connected, we may therefore choose maximal compact subgroups
$C_K,C_L\leq H$ such that
\[
K\subseteq C_K,
\qquad
L\subseteq C_L.
\]
There is then an element $a\in H$ such that
\begin{equation}\label{eq:conjugate_subgroups}
aC_Ka^{-1}=C_L.
\end{equation}
The groups $C_K$ and $C_L$ are compact and hence closed in $H$, so
the closed subgroup theorem also makes them compact Lie groups.

The inclusion $H\to G$ is continuous. Consequently, the images of
$C_K$ and $C_L$ are compact subgroups of
$\operatorname{Aut}(D)$. Every compact group of biholomorphic
automorphisms of a bounded convex domain has a common fixed point; see
\cite[Theorem~4.1.2, pp.~101--102]{GreeneKimKrantz2011}.

It follows that $C_K$ has a common fixed point $x\in D$. Since
$K\subseteq C_K$, the point $x$ is fixed by $K$. The hypothesis
$\operatorname{Fix}(K)=\{0\}$ therefore implies $x=0$. Thus
$C_K$ fixes $0$. Similarly, $C_L$ fixes $p$.

We now use \eqref{eq:conjugate_subgroups}. If $c\in C_L$, write $c=ac_Ka^{-1}$ for some
$c_K\in C_K$. Since $C_K$ fixes $0$,
\[
c(a(0))
=
ac_Ka^{-1}(a(0))
=
ac_K(0)
=
a(0).
\]
Hence $a(0)$ is fixed by every element of $C_L$, and in particular
by every element of $L$. Since
$\operatorname{Fix}(L)=\{p\}$, we obtain
\[
a(0)=p.
\]
Therefore
\[
a^{-1}(p)=0.
\]

Finally, $a^{-1}\in H=\langle K,L\rangle$. By the definition of the
generated subgroup, $a^{-1}$ is a finite word in elements of $K$
and $L$. After multiplying consecutive factors from the same subgroup
and inserting identity factors if necessary, we may write
\[
a^{-1}=k_1\ell_1\cdots k_m\ell_m,
\qquad
k_j\in K,\quad \ell_j\in L.
\]
Taking $w=a^{-1}$ proves the result.
\end{proof}

\subsection{From isomorphisms to vacuum preserving isomorphisms}

Our goal is now, given an isomorphism $\alpha \colon A(\fV) \to A(\fW)$, to construct a vacuum preserving isomorphism $\beta \colon A(\fV) \to A(\fW)$. 
It will be convenient to use the notion of the {\em singular nucleus} of a homogeneous variety from \cite{DRS1}, defined as follows. 

If $V$ is a homogeneous variety in a ball (meaning that $V = B \cap Z(I_V)$, where $B$ is the unit ball with respect to some norm and $Z(I_V)$ is the zero locus of some homogeneous ideal $I_V \triangleleft \C[z_1,\ldots,z_d]$) then the singular locus $S(V)$ of $V$ is non-empty unless $V$ is a linear subspace. 
When $S(V)$ is non-empty, it is a homogeneous variety in its own right, and therefore has its own singular locus. 
The {\em singular nucleus} of $V$ is the linear subspace $N = S(\cdots S(V) \cdots)$ where the singular locus operator $S(\cdot)$ is applied until one gets a linear subspace (if $V$ is already a linear subspace then define the singular nucleus $N$ to be $V$). 
It is clear that a biholomorphic map between two homogeneous varieties must preserve its singular nucleus. 

We continue the standing assumption that $E$ and $F$ are two finite-dimensional operator space structures. 
Note that if $\fV\subseteq \bB_E$ is a homogeneous NC variety in $\bB_E$ then
$V:=\fV_1$ is a homogeneous variety in the ball $\bB_E(1)$. 

\begin{theorem}\label{thm:iso_vacuum_iso}
Let $\fV \subseteq \bB_E$ and $\fW \subseteq \bB_F$ be two homogeneous varieties, and let $\alpha \colon A(\fV) \to A(\fW)$ be a bounded isomorphism. 
Then there exists a bounded isomorphism $\beta \colon A(\fV) \to A(\fW)$ such that $\beta^*\Phi_0 = \Phi_0$. 
If, in addition, $\alpha$ 
and $\alpha^{-1}$ are
completely bounded, then $\beta$ can be chosen completely bounded.
\end{theorem}
\begin{proof}
    Given a (completely) bounded isomorphism $\alpha \colon A(\fV) \to A(\fW)$, we know by Proposition \ref{prop:homo_ind_NCmap} that $\alpha^*$ induces an NC biholomorphism $G = \alpha^* \big|_{\widetilde{\fW}} \colon \widetilde{\fW} \to \widetilde{\fV}$. 
    The map $G$ restricts to a biholomorphism $g \colon W = \fW_1 \to V = \fV_1$ between the first level scalar varieties. 
    Let $N$ be the singular nucleus of $W$ and $M$ the singular nucleus of $V$. Then since $g$ is a biholomorphism, $g$ maps $N$ onto $M$. Now, $M,N$ are intersections of linear subspaces with unit balls, hence $M$ and $N$ can both be considered as bounded convex open and balanced domains in some $\bC^k$. 

    Now let $\gamma_t \colon X \mapsto e^{it} X$ and $\lambda_t : X \mapsto e^{it} X$ be the circle action on $\fV$ and $\fW$, respectively. Note that $\gamma_t = \Gamma_t^*\big|_{\fV}$ is the map induced by the completely isometric automorphism $\Gamma_t f (X) = f(e^{itX})$ on $A(\fV)$ and likewise $\lambda_t$ is induced by a rotation automorphism $\Lambda_t$ on $A(\fW)$. 

    Suppose that $G(0_1)=g(0) = p \neq 0$, for otherwise we are done. 
    Consider the compact connected groups $K$ and $L$ generated in $\operatorname{Aut}(M)$ by (the restriction of) $\gamma_t$ and $g \circ \lambda_t \circ g^{-1}$. 
    By definition, the fixed point of $K$ is $0$ and the fixed point of $L$ is $p$. 
    By Lemma \ref{lem:two_circ} there is a finite word of the form 
    \[
    w = \gamma_{s_1} \circ g \circ \lambda_{t_1} \circ g^{-1} \circ \cdots \circ \gamma_{s_m} \circ g \circ \lambda_{t_m} \circ g^{-1}
    \]
    such that $w(p) = 0$. 
    But
    \[
    w = \Gamma_{s_1}^* \circ \alpha^* \circ \Lambda^*_{t_1} \circ \cdots \Lambda_{t_m}^* \circ (\alpha^{-1})^*  
    = \left( \alpha^{-1} \circ \Lambda_{t_m} \circ \cdots \circ \alpha \circ \Gamma_{s_1}\right)^* .
    \]
    Now define 
    \[
    \varphi = \alpha^{-1} \circ \Lambda_{t_m} \circ \cdots \circ \alpha \circ \Gamma_{s_1}, 
    \]
    so that $w = \varphi^*$. 
    It follows from the definitions that 
    \[
    \beta := \alpha \circ \varphi
    \]
    is a bounded isomorphism from $A(\fV)$ onto $A(\fW)$ and (identifying points with their corresponding evaluation representations)
    \[
    \beta^* (0) = \varphi^* \circ \alpha^* (0) = w \circ G(0) = w(p) = 0, 
    \]
    so that $\beta^* \Phi_0 = \Phi_0$. Moreover, if $\alpha$ and $\alpha^{-1}$ are completely bounded then so are $\beta$ and $\beta^{-1}$.

\end{proof}

Thanks to the previous theorem, in order to show that the existence of an isomorphism $\alpha \colon A(\fV) \to A(\fW)$ implies the existence of a linear map $\widetilde{\fW} \to \widetilde{\fV}$, it suffices to do this under the assumption that $\alpha$ is vacuum preserving. 
We first show that the existence of a vacuum preserving isomorphism implies the existence of an isomorphism $\beta \colon A(\fV) \to A(\fW)$ which is {\em graded}, in the sense that $\beta$ maps homogeneous polynomials of degree $n$ to homogeneous polynomials of degree $n$. 
The idea behind the following proof is related to the maneuver used in \cite{dor2014operator} to pass from ``semi-graded" to graded isomorphisms in the theory of subproduct systems \cite{SS09} (see also \cite{SalinasPillajo}). 

\begin{theorem}\label{thm:graded}
Let $\fV \subseteq \bB_E$ and $\fW \subseteq \bB_F$ be two homogeneous varieties, and let $\alpha \colon A(\fV) \to A(\fW)$ be a bounded isomorphism which is vacuum preserving, in the sense that $\alpha^*\Phi_0 = \Phi_0$. 
Then there exists a graded isomorphism $\beta \colon A(\fV) \to A(\fW)$ such that $\|\beta\|_\kappa\leq\|\alpha\|_\kappa, \, \|\beta^{-1}\|_\kappa\leq\|\alpha^{-1}\|_\kappa$ for $\kappa \in \{b,cb\}$.  
\end{theorem}

\begin{proof}
Let $P_n$ be the projection \eqref{eq:homog_proj} onto the $n$th homogeneous component. Define
$
\cE_n:=P_nA(\fV)
$
and
$
\cF_n:=P_nA(\fW).
$
For every $n\geq0$ and $f_n\in\cE_n$, define
$\alpha_n:\cE_n\longrightarrow\cF_n$ by
\[
\alpha_n(f_n)
:=P_n\alpha(f_n)
=
\frac{1}{2\pi}\int_0^{2\pi}
e^{-int}\Lambda_t\big(\alpha(f_n)\big)\,dt.
\]
Define $\beta$ on finite sums of homogeneous elements by
$$
\beta\left(\sum_{n=0}^N f_n\right)
:=
\sum_{n=0}^N\alpha_n(f_n)
=
\frac{1}{2\pi}\int_0^{2\pi}
\Lambda_{-t}\alpha\Gamma_t
\left(\sum_{n=0}^N f_n\right)\,dt,
\qquad f_n\in\cE_n,
$$
where $\Gamma_t$ and $\Lambda_t$ denote the gauge automorphisms of $A(\fV)$ and $A(\fW)$, respectively.
Since $\Gamma_t$ and $\Lambda_t$ are completely isometric automorphisms, for every finite sum $p=\sum_{n=0}^N f_n$, we have
\begin{align*}
\|\beta(p)\|
&\leq \frac{1}{2\pi}\int_0^{2\pi}
\|\Lambda_{-t}\alpha\Gamma_t(p)\|\,dt
\leq \|\alpha\|\,\|p\|.
\end{align*}
Now, for $f\in A(\fV)$, the Cesàro sums
$\sigma_N(f)\to f$ in $A(\fV)$. Moreover,
$
\|\beta(\sigma_N(f))-\beta(\sigma_M(f))\|
\leq
\|\alpha\|\,\|\sigma_N(f)-\sigma_M(f)\|.
$
Thus $(\beta(\sigma_N(f)))_{N\geq0}$ is a Cauchy sequence in $A(\fW)$. We may therefore extend $\beta$ to $A(\fV)$ by defining
$
\beta(f):=\lim_{N\to\infty}\beta(\sigma_N(f)).
$
By continuity and the integral representation on finite homogeneous sums, we obtain
$$
\beta(f)
=
\frac{1}{2\pi}\int_0^{2\pi}
\Lambda_{-t}\alpha\Gamma_t(f)\,dt,
\qquad f\in A(\fV).
$$
Consequently, $\beta$ is a bounded linear map with
$
\|\beta\|\leq\|\alpha\|.
$
We next show that $\beta$ is multiplicative. Define
$$
\cI_{\fV}:=\ker\Phi_0\subseteq A(\fV),
\qquad
\cI_{\fW}:=\ker\Phi_0\subseteq A(\fW),
$$
where $\Phi_0$ denotes evaluation at the origin. For every $n\geq1$, let
$$
\cI_{\fV}^n
:=
\spn\{f_1\cdots f_n:
f_j\in\cI_{\fV},\ 1\leq j\leq n\},
$$
and denote its norm closure by
$\overline{\cI_{\fV}^n}$.
Define $\cI_{\fW}^n$ and
$\overline{\cI_{\fW}^n}$ similarly.
Since $\alpha$ is a vacuum-preserving isomorphism,
$
\alpha(\cI_{\fV})=\cI_{\fW}.
$
By multiplicativity and continuity of $\alpha$, we obtain
$$
\alpha\big(\overline{\cI_{\fV}^n}\big)
=
\overline{\cI_{\fW}^n},
\qquad n\geq1.
$$
Moreover, for every $n\geq1$,
$$
\overline{\cI_{\fV}^n}
=
\{f\in A(\fV):P_kf=0,\ 0\leq k<n\},
$$
and the analogous identity holds for $\fW$.
Indeed, every product of $n$ functions vanishing at the origin has no homogeneous components of degree less than $n$. Conversely, the Cesàro sums of a function whose first $n$ homogeneous components vanish consist of homogeneous polynomials of degree at least $n$, each of which belongs to $\cI_{\fV}^n$. The identity follows by norm convergence of the Cesàro sums.
Consequently, $\alpha$ and $\alpha^{-1}$ preserve the homogeneous-degree filtrations. In particular, for every $f_n\in\cE_n$,
$$
P_k\alpha(f_n)=0,\qquad k<n.
$$
We next show that the maps $\alpha_n$ respect multiplication. Let
$n,m\geq1$, $f\in\cE_n$, and $h\in\cE_m$. By the preceding discussion,
$$
\alpha(f)=\alpha_n(f)+f'
\quad
(f'\in\ol{\cI_{\fW}^{(n+1)}}),
\qquad
\alpha(h)=\alpha_m(h)+h'
\quad
(h'\in\ol{\cI_{\fW}^{(m+1)}}).
$$
Since $\alpha$ is multiplicative,
$
\alpha(fh)=\alpha(f)\alpha(h).
$
Taking the homogeneous component of degree $n+m$, all terms involving
$f'$ or $h'$ vanish, and hence
$
\alpha_{n+m}(fh)=\alpha_n(f)\alpha_m(h).
$
If $n=0$ or $m=0$, the same identity follows immediately from
$\alpha_0=I_{\mathbb C1}$.
Hence $\beta$ is multiplicative on finite sums of homogeneous elements, and therefore on $A(\fV)$ by continuity and Cesàro approximation.

Applying the same construction to $\alpha^{-1}$ gives a bounded graded homomorphism
$\delta:A(\fW)\longrightarrow A(\fV)$
such that
$$
\delta|_{\cF_n}
=
P_n\alpha^{-1}|_{\cF_n},
\qquad
\|\delta\|\leq\|\alpha^{-1}\|.
$$
For every $f_n\in\cE_n$, since
$
\alpha(f_n)-P_n\alpha(f_n)
\in\overline{\cI_{\fW}^{n+1}}
$
and $\alpha^{-1}$ preserves the homogeneous-degree filtration, we have
\[
\begin{aligned}
\delta\beta(f_n)
=P_n\alpha^{-1}\big(P_n\alpha(f_n)\big)=P_n\alpha^{-1}\alpha(f_n)
=f_n.
\end{aligned}
\]
Similarly, $\beta\delta(g_n)=g_n$ for every
$g_n\in\cF_n$. By linearity and density of homogeneous polynomials,
$\delta\circ\beta=I$ and
$\beta\circ\delta=I$
on the respective algebras.
Thus $\delta=\beta^{-1}$, and hence
$
\|\beta^{-1}\|\leq\|\alpha^{-1}\|.
$
Therefore $\beta$ is a bounded graded isomorphism.

Finally, if both $\alpha$ and $\alpha^{-1}$ are completely bounded, the same averaging argument applied at every matrix level gives
$
\|\beta\|_{cb}\leq\|\alpha\|_{cb},
$
and
$
\|\beta^{-1}\|_{cb}\leq\|\alpha^{-1}\|_{cb}.
$
This completes the proof.
\end{proof}

We come to the following sharpened classification theorem. 
We remind the reader that $A \colon \bC^e \to \bC^d$ is a linear map between finite-dimensional operator spaces, then $A$ can be promoted to an NC map $A \colon \bM(\bC^e) \to \bM(\bC^d)$ by acting as $I_n \otimes A$ on $M_n(\bC^e) = M_n \otimes \bC^e$, that is 
\[
A (X_1, \ldots, X_e) = \left(\sum_{j=1}^e A_{1j} X_j, \ldots, \sum_{j=1}^e A_{dj} X_j \right). 
\]

\begin{theorem}\label{thm:classification_linear}
Let $\fV\subseteq\bB_E$ and $\fW\subseteq\bB_F$
be homogeneous NC varieties, and let
$\kappa\in\{b,cb\}$.
Then $A(\fV)$ and $A(\fW)$ are
$\kappa$-boundedly graded isomorphic if and only if
there exist linear maps $L:F\to E$ and $M:E\to F$
whose matrix amplifications induce mutually inverse
$\delta_\kappa$-bi-Lipschitz NC linear maps
$L:\widetilde{\fW}\longrightarrow\widetilde{\fV}$ and $
M:\widetilde{\fV}\longrightarrow\widetilde{\fW}$.
\end{theorem}

\begin{remark}
Such linear maps $L$ and $M$ are necessarily completely bounded hence uniformly continuous, and so once they are assumed to be $\delta_\kappa$-bi-Lipschitz it follows that they induce a $\kappa$-bounded graded isomorphism by Theorem \ref{thm:bounded-cb-classification}.  
\end{remark}
\begin{proof}
By the above remark, we need only prove that the existence of a $\kappa$-bounded isomorphism implies the existence of such $\delta_\kappa$-bi-Lipschitz linear maps. 
By the combination of Theorems \ref{thm:iso_vacuum_iso} and \ref{thm:graded}, we may as well assume that $\alpha$ is graded. 
This means that $h_j = \alpha(Z_j)$ is a homogeneous polynomial of degree $1$ on $\fW$ for all $j=1, \ldots, d$. 
We may think of each $h_j$ as a polynomial function of degree $1$ on the space. 
Let $L$ be the linear map given by 
\[
Lw = (h_1(w), \ldots, h_d(w)). 
\]
We know by Theorem \ref{thm:bounded-cb-classification} that there exists a bi-Lipschitz map $G \colon \widetilde{W} \to \widetilde{\fV}$ such that 
\[
\alpha(f) = \widetilde{f}\circ G, 
\]
for all $f \in A(\fV)$. 
We shall show that $G = L\big|_{\widetilde{\fW}}$. 
But for every $j = 1, \ldots, d$, the $j$th coordinate of $G$ is given by 
\[
G_j = \widetilde{Z_j} \circ G = \alpha(Z_j) = h_j, 
\]
and this shows that $L$ gives a bijection of $\widetilde{\fW}$ onto $\widetilde{\fV}$. 
Similarly, one obtains $M$ from the graded isomorphism $\alpha^{-1}$. 
The fact that $L$ and $M$ restrict to mutual inverses follows from functoriality. 
\end{proof}

\begin{remark}\label{rem:rmk_homo}
In Theorem 10.3 in \cite{SSS2} a classification of NC function algebras on homogeneous subvarieties of the row ball $\fB_d$ was obtained, which appears somewhat sharper than the above theorem.
It was shown in that restricted setting that, as above, the existence of a $\kappa$-bounded isomorphism is equivalent to the existence of $\delta_\kappa$-bi-Lipschitz map; however in \cite{SSS2} it was also shown that a bounded graded isomorphism is automatically completely bounded. Two natural questions come to mind. 

First, is it true in our general setting that a bounded isomorphism must be completely bounded. The answer is no: in Section \ref{subsec:bdd_not_cb} below we provide an example of to isometrically isomorphic NC function algebras on homogeneous varieties that are not completely boundedly isomorphic. 

Second, linear maps between finite-dimensional spaces are such well-behaved objects, so one wonders whether one really requires an assumption that a linear bijection between similarity envelopes of homogeneous varieties be bi-Lipschitz. 
In fact, this was an important problem left open in \cite{SSS2}. 
It was tempting to conjecture that linear bijections between similarity envelopes of homogeneous varieties are automatically bi-Lipschitz, because that generalizes nicely the classification of complete Pick algebras in the homogeneous setting (see \cite{salomon2016isomorphism} for an overview). 
Indeed, if $\cM_d$ denotes the multiplier algebra of the Drury--Arveson space and $\cA_V$ is the closure of polynomials in $\cM_d\big|_V$, then it was shown in \cite{Hartz} (following \cite{DRS1}) that for two homogeneous varieties $V,W \subseteq \bB_d$ $\cA_V$ and $\cA_W$ are isomorphic if and only if they are completely boundedly isomorphic, and that this happens if and only if there exists a linear map which is a bijection $W \to V$. 
Now, by Sections 9 and 10 in \cite{SSS2} the algebras $\cA_V$ can be identified with $A(\fV)$ for an appropriate reduced affiliated NC subvariety $\fV$ of the row ball $\fB_d$. 
It follows by a combination of Hartz's result and Theorem \ref{thm:classification_linear} that a linear map between two homogeneous subvarieties of the Euclidean ball extends to a bi-Lipschitz mapping between the similarity envelopes of the reduced affiliated NC varieties. 
So it is natural to conjecture that a linear bijection between similarity envelopes is automatically bi-Lipschitz. 
This turns out to be false, as we show in Section \ref{subsec:cnter_examples2} below. 
\end{remark}

\subsection{Linearization of maps}

The following theorem should be compared with Theorem \ref{thm:bounded-cb-classification}. 
It shows, that already the existence of a bi-Lipschitz NC bijection between similarity envelopes of homogeneous varieties implies that the algebras are (completely) boundedly isomorphic. 

\begin{theorem}
\label{thm:linearisation_maps}
Let $\fV\subseteq\bB_E$ and $\fW\subseteq\bB_F$ be
homogeneous NC varieties. Then
$A(\fV)$ is (completely) boundedly isomorphic to $A(\fW)$ if and
only if there exists an NC bijection
$G:\widetilde{\fW}\longrightarrow\widetilde{\fV}$
which is ($\delta_{cb}$-bi-Lipschitz) $\delta_b$-bi-Lipschitz.
\end{theorem}

\begin{proof}
One direction follows from Theorem \ref{thm:bounded-cb-classification}. For the converse,
fix $\kappa\in\{b,cb\}$ and suppose that
$G:\widetilde{\fW}\to\widetilde{\fV}$ is an NC
$\delta_\kappa$-bi-Lipschitz bijection, with inverse
$H=G^{-1}$.

As in the proof of Theorem \ref{thm:iso_vacuum_iso}, we may use the two-circle trick (Lemma \ref{lem:two_circ}) to replace 
$G$ by an NC bijection
$F_0:\widetilde{\fW}\to\widetilde{\fV}$ satisfying
$F_0(0)=0$. Let $H_0=F_0^{-1}$. Then $F_0$ and $H_0$ are mutually inverse
$\delta_\kappa$-bi-Lipschitz NC maps satisfying
$F_0(0_n)=H_0(0_n)=0_n$ for $n\geq1.$

By the proof of the implication $(1) \Rightarrow (2)$ in Proposition~\ref{prop:equiv_G_norm}, we have that 
$
C_{F_0}^{\kappa}
:=
\sup_{X\in\fW}\|\Phi_{F_0(X)}\|_\kappa
<\infty$, and 
$C_{H_0}^{\kappa}
:=
\sup_{Y\in\fV}\|\Phi_{H_0(Y)}\|_\kappa<\infty
$. 
Consequently, we get a bounded unital
homomorphism $T_{F_0}:A(\fV)\longrightarrow H^\infty(\fW)$ given by
$$T_{F_0}(f)=\widetilde f\circ(F_0|_{\fW}), \qquad f \in A(\fV) $$
with
$\|T_{F_0}\|_\kappa\leq C_{F_0}^{\kappa}.$ Similarly, composition with $H_0$ defines a bounded
homomorphism
$T_{H_0}:A(\fW)\longrightarrow H^\infty(\fV)$ by $
T_{H_0}(h)
=
\widetilde h\circ (H_0|_{\fV})$ for $ h \in A(\fW) $,
with
$\|T_{H_0}\|_\kappa\leq C_{H_0}^{\kappa}.$

Since $F_0|_{\fW}$ and $H_0|_{\fV}$ have
bounded NC coordinate functions and the varieties
are homogeneous, their degree-one homogeneous
components are represented by ambient linear maps
$A:F\to E$ and $B:E\to F$, respectively.
As $F_0(0)=H_0(0)=0$, their homogeneous expansions
take the form
\[
F_0(X)=AX+\sum_{k\geq2}(F_0)_k(X),
\qquad
H_0(Y)=BY+\sum_{k\geq2}(H_0)_k(Y).
\]

Fix $X\in\fW$ and $Y\in\fV$. Since $F_0$ and $H_0$
are continuous at zero and
$F_0(0)=H_0(0)=0,$
for all sufficiently small $t$ we have
$H_0(tY) \in \widetilde{\fW}\cap\bB_F=\fW$ and $F_0(tX)\in \widetilde{\fV}\cap\bB_E= \fV.$
Since
$F_0(H_0(tY))=tY,
\,
H_0(F_0(tX))=tX$,
taking first-order terms yields
\[
ABY=Y ,
\qquad
BAX=X .
\]

For a polynomial function $p\in\operatorname{Pol}(\fV)$,
let $p=\sum_{n=0}^N p_n$ be its homogeneous
decomposition, and define
\[
\beta(p)=
\sum_{n=0}^N P_n\bigl(T_{F_0}(p_n)\bigr).
\]
If $q=\sum_{n=0}^N q_n$ is an ambient polynomial
representative of $p$, then, since $F_0(0)=0$, we obtain $P_n\bigl(T_{F_0}(p_n)\bigr)(X)=q_n(AX).$
Consequently,
\[
\beta(p)(X)=q(AX),\qquad X\in\fW.
\]
The homogeneous decomposition is unique on
$\fV$, so $\beta$ is well defined. The above
formula also shows that $\beta$ is a graded
homomorphism with polynomial range.
Let $\Gamma_t$ and $\Lambda_t$ denote the
gauge actions on $A(\fV)$ and $H^\infty(\fW)$,
respectively. Then
\[
\beta(p)=\frac{1}{2\pi}
\int_0^{2\pi}
\Lambda_{-t}T_{F_0}\Gamma_t(p)\,dt.
\]
Since the gauge actions are completely isometric,
\[
\|\beta(p)\|_\kappa
\leq \|T_{F_0}\|_\kappa\|p\|_\kappa
\leq C_{F_0}^{\kappa}\|p\|_\kappa,
\]
where the estimate is understood at all matrix
levels when $\kappa=cb$.
Thus $\beta$ extends by polynomial density to a
$\kappa$-bounded graded homomorphism
$\beta:A(\fV)\longrightarrow A(\fW).$

Similarly, $H_0$ yields a $\kappa$-bounded
graded homomorphism
$\eta:A(\fW)\to A(\fV)$,
which acts on polynomial functions by
composition with $B$.
Since $ABY=Y$ for $Y\in\fV$ and
$BAX=X$ for $X\in\fW$, we have
\[
\eta\beta(p)=p,\qquad
\beta\eta(r)=r
\]
for all $p\in\operatorname{Pol}(\fV)$ and
$r\in\operatorname{Pol}(\fW)$.
By polynomial density,
\[
\eta\beta=I_{A(\fV)},
\qquad
\beta\eta=I_{A(\fW)}.
\]
Hence $\beta$ is a $\kappa$-bounded algebra
isomorphism.
\end{proof}

\section{Examples}\label{sec:examples}

In this section, we collect examples demonstrating the necessity of various assumptions made in the results of the previous section.

\subsection{Lipschitz versus uniformly continuous}\label{subsec:cnter_examples1}

The following two examples show that ambient uniform continuity and $\delta_\kappa$-Lipschitz continuity are, in general, independent conditions. 
Thus it is really necessary to assume both in the classification Proposition \ref{prop:equiv_G_norm}.

\begin{example}
 Let
$\fV=\fW=\bB_{\C}.$
At the scalar level this is the unit disc $\bD$, and
$\widetilde{\bB_{\C}}(n)
=
\{X\in M_n(\C):\rho(X)<1\}.$ 
Consider
\[
g(z)
=
\exp\left(-\frac{1+z}{1-z}\right),
\qquad z\in\bD.
\]
Since
$\operatorname{Re}\frac{1+z}{1-z}>0
$ for $z\in\bD,$
we have $|g(z)|<1$.
Hence, by holomorphic functional calculus,
$G(X):=g(X)$
defines an NC holomorphic self-map of
$\widetilde{\bB_{\C}}$.
For every
$F\in H^\infty(\bD)_1$,
$F\circ g\in H^\infty(\bD)_1.$ 
Hence, Proposition~\ref{prop:equi_F_norm} yields
$\delta_k(G(X),G(Y))
=\sup_{F\in H^\infty(\bD)_1}
\|\widetilde{F\circ g}(X)
-\widetilde{F\circ g}(Y)\| \leq \delta_k(X,Y),$
for $\kappa\in\{b,cb\}$.
Thus $G$ is $\delta_\kappa$-Lipschitz.
On the other hand, $g$ is a well-known example of an inner function which is not uniformly continuous, even at the scalar level. 
\end{example}

\begin{example}\label{ex:unif_not_bilip}
Conversely, a map $G \colon \widetilde{{\fW}} \to \widetilde{{\fV}}$ with components in $A(\fW)$ need not be $\delta_b$-Lipschitz. 
Indeed, let $\fV = \bB_E$ be the $d$-dimensional row ball and let $\fW = \bB_F$ be the $d$-dimensional column ball. Consider  $G: \widetilde{\bB}_F \to \widetilde{\bB}_E$ given by $G(X) = X$.
It follows from \cite[Example 5.2]{SchShaSha} that the identity map $G(X) = X$ is a bijection $\widetilde{\bB}_F \to \widetilde{\bB}_E$. 
Since the identity gives a completely bounded map between the row and column operator spaces, $G$ restricts to a uniformly continuous map on $\bB_F$ (alternatively, as
$G_j|_{\bB_F}=Z_j\in A(\bB_F),
$ for $ 1 \leq j \leq d,$ Lemma \ref{lem:uniform_components} implies that
$G|_{\bB_F}$ is uniformly continuous).  
However, by \cite[Remark 8.5]{SchShaSha}, the induced mapping $p \mapsto p \circ G$ on polynomials is not bounded. By Proposition \ref{prop:equiv_G_norm}, $G$ is not $\delta_b$-Lipschitz nor is it $\delta_{cb}$-Lipschitz. 
\end{example}

\subsection{Linear bijections between varieties are not automatically bi-Lipschitz}\label{subsec:cnter_examples2}

By Theorem \ref{thm:classification_linear}, two algebras of uniformly continuous NC functions on homogeneous varieties are (completely) boundedly isomorphic if and only if there is a linear map that restricts to a bijection between the similarity envelopes of the varieties, which is bi-Lipschitz with respect to the appropriate pseudohyperbolic metric. 
Linear maps between finite dimensional spaces are automatically completely bounded, so it is natural to ask whether the bi-Lipschitz condition is automatic as well. 
Such hopes are grounded on experience: by work on the isomorphism problem for complete Pick multiplier algebras on homogeneous varieties \cite{DRS1,Hartz,Hartz2} the automatic-bi-Lipschitzness property holds in the case of radical commutative varieties in the row ball (see also \cite[Remark 9.6 and Example 9.9]{SSS2}). 
The following example shows that this property fails for non-radical commutative subvarieties of the row ball.

\begin{example}\label{ex:cntr_exp_BLip}
    In this example we consider commuting varieties cut out in the row ball. 
    We shall exhibit a linear bijection $A$ between the similarity envelopes of two varieties in the row ball which is not bi-Lipschitz (thereby answering the second question raised in Remark \ref{rem:rmk_homo}). In fact, the $A$ we exhibit will be spectral radius preserving on tuples. 
    
    Let us consider the free algebra $\C \langle x, y \rangle$. 
    Consider the ideals
    \[
    I = \langle xy - yx, (x-y)^2 \rangle \subset \C \langle x, y \rangle, 
    \]
    and 
    \[
    J = \langle xy - yx, (x-\sqrt{3}y)^2 \rangle \subset \C \langle x, y \rangle.
    \]
    The quotient by $I$ corresponds to a double line $x = y$. 
    Let the variety cut out in the row ball by $I$ be 
    \[
    \fV = V_{\fB_2}(I) := \{X \in \fB_2 : p(X) = 0 \textrm{ for all } p \in I\}. 
    \]
   Similarly let the variety cut out by $J$ be $\fW = V_{\fB_2}(J)$. 

    Define the matrix
    \[
    A = \frac{1}{\sqrt{2}} \begin{pmatrix} \sqrt{3} & 0 \\ 0 & 1 \end{pmatrix} .
    \]
    Now note that $(z,z) \in \bB_2$ if and only if $|z| < \frac{1}{\sqrt{2}}$, and this happens if and only if $\frac{1}{\sqrt{2}} (\sqrt{3} z, z) \in \bB_2$. Hence, $A$ maps $\fV(1)$ onto $\fW(1)$ isometrically. 
    We claim that $A$ maps $\widetilde{\fV}$ onto $\widetilde{\fW}$. 
    Indeed, $A$ maps the NC zero locus $Z_{nc}(I)$ of $I$ bijectively onto the NC zero locus $Z_{nc}(J)$ of $J$. 
    Moreover, if $(X,Y) \in \widetilde{\fV}$, then $X$ and $Y$ commute. 
    So up to a unitary, we may assume that $X$ and $Y$ are upper triangular 
    and the diagonal points are in $\fV(1)$. Hence, $A(X,Y)$ has the same spectral radius as $(X,Y)$. Using that $\widetilde{\fB}_2 = \{(X,Y) : \rho_{\operatorname{row}}(X,Y)<1\}$, it follows that $(X,Y) \in \widetilde{\fV}$ if and only if $A(X,Y) \in \widetilde{\fW}$. 

    For every $k \in \mathbb N$, we define a point $(X_k,Y_k)\in\widetilde {\fV}(2)$ as follows:
    \[
    X_k = \frac{1}{\sqrt{2k^2 +1}} \begin{pmatrix} k & 1/2 \\ 0 & k \end{pmatrix} \text{ and } Y_k = \frac{1}{\sqrt{2k^2 + 1}} \begin{pmatrix} k & - 1/2 \\ 0 & k \end{pmatrix}.
    \]
    Let us compute the norm of the point
    \begin{align*}
    X_k X_k^* + Y_k Y_k^* 
    &= 
    \frac{1}{2k^2 + 1} \left( \begin{pmatrix} k^2 + 1/4 & k/2 \\ k/2 & k^2 \end{pmatrix} + \begin{pmatrix} k^2  + 1/4 & -k/2 \\ -k/2 & k^2 \end{pmatrix} \right) \\ 
    &= \frac{1}{2k^2 + 1} \begin{pmatrix} 2 k^2 + 1/2 & 0 \\ 0 & 2k^2\end{pmatrix}.
    \end{align*}
    Hence, 
    \[
    \|(X_k, Y_k)\|_{row} = \sqrt{\frac{2 k^2 + 1/2}{2 k^2 + 1}} < 1. 
    \]
    It is also clear that $(X_k, Y_k)$ belong to the variety on the second level, since $X_k - Y_k$ is square nilpotent.
    However, take the linear polynomial
    \[
    \ell(x,y) = \frac{1}{2} (\sqrt{3} x + y).
    \]
    Then, $\|\ell\|_{\infty} = 1$. However,
    \[
    \ell(A(X_k,Y_k)) = \frac{1}{2  \sqrt{2} \sqrt{2 k^2+1}} \begin{pmatrix} 4 k & 1 \\ 0 & 4 k \end{pmatrix} = \underbrace{\frac{\sqrt{2} k}{\sqrt{2 k^2 + 1}}}_{t_k} I + \frac{1}{2  \sqrt{2} \sqrt{2 k^2+1}} \underbrace{\begin{pmatrix} 0 & 1 \\ 0 & 0 \end{pmatrix}}_{N}.
    \]
    Let $\varphi_k(z) = \frac{z - t_k}{1 - t_k z}$. Note that since $t_k < 1$, $\varphi_k$ is an automorphism of the unit disc. Since $\ell$ is a norm 1 homogeneous polynomial, it is an isometry on the Fock space. Hence, so is $\varphi_k(\ell)$. In particular, $\|\varphi_k(\ell)\|_{\infty} = 1$. Now,
    \begin{multline*}
    \|\varphi_k(\ell(A(X_k,Y_k)))\| = \frac{1}{2  \sqrt{2} \sqrt{2 k^2+1}} \left\| N \left(I - t_k \ell(A(X_k,Y_k))\right)^{-1} \right\| =  \\ \frac{1}{2  \sqrt{2} \sqrt{2 k^2+1}} \left\|N \left( (1 - t_k^2) I - \frac{t_k}{2  \sqrt{2} \sqrt{2 k^2+1}} N \right)^{-1}\right\| = \frac{1}{2  \sqrt{2} \sqrt{2 k^2+1}} \frac{1}{1 - t_k^2}.
    \end{multline*}
    Now $t_k^2 = \frac{2 k^2}{2 k^2 + 1} = 1 - \frac{1}{2 k^2 + 1}$ and so $1 - t_k^2 = \frac{1}{2 k^2 + 1}$. Hence,
    \[
    \lim_{k \to \infty} \|\varphi_k(\ell(A(X_k,Y_k)))\| = \lim_{k \to \infty} \frac{\sqrt{2k^2 + 1}}{2  \sqrt{2}} = \infty.
    \]
    Thus, Condition (2) of Proposition \ref{prop:equiv_G_norm} is not satisfied for $\kappa = b$, and it follows that $A$ is not bi-Lipschitz with respect to either metric.
\end{example}

\subsection{Bounded isomorphisms need not be completely bounded}\label{subsec:bdd_not_cb}

In Remark \ref{rem:rmk_homo}, we asked whether bounded isomorphisms between NC function algebras on homogeneous varieties must be completely bounded, too. 
In fact, one might also ask a somewhat different question, whether being boundedly isomorphic is equivalent to being completely boundedly isomorphic. 
The following example answers both question in the negative, by showing that there exist isometric NC function algebras on homogeneous varieties that are are not completely boundedly isomorphic. 
This is in contrast with \cite[Theorem 10.3]{SSS2} where it was shown that if we know that both varieties are in the row ball then bounded isomorphism implies completely bounded.

\begin{example}\label{ex:bdd_not_cb}
Let $E = (\ell^2_d)_{\operatorname{row}}$ and $F = (\ell^2_d)_{\operatorname{col}}$ be the row and column operator space structures on $\bC^d$, with $d \geq 2$. 
Let $\fC \bB_E$ and $\fC \bB_F$ be the {\em commutative NC balls}, that is 
\[
\fC \bB_E = \{X = (X_1, \ldots, X_d) \in \bB_E : X_i X_j = X_j X_i \textrm{ for all } i,j=1,\ldots, d\}, 
\]
and similarly for $\fC\bB_F$. 
We claim that $\fV = \fC \bB_E$ and $\fW = \fC \bB_F$ are homogeneous varieties for which the natural identity map $A(\fV) \to A(\fW)$ is an isometric algebra isomorphism but is not completely bounded.

Indeed, $X \in \fC\bB_E$ if and only if the transpose
\[
X^T := (X_1^T, \ldots, X_d^T) \in \fC\bB_F, 
\]
and so using the fact that the transpose is isometric on matrices and preserves commutativity, we have for every polynomial $p \in \bC\langle Z_1, \ldots, Z_d \rangle$
\begin{align*}
    \|p\|_{A(\fC\bB_E)} 
    &= \sup\{\|p(X)\| : X \in  \fC\bB_E\ \} 
    = \sup\{\|p(X)^T\| : X \in  \fC\bB_E\ \} \\
    &= \sup\{\|p(X^T)\| : X \in  \fC\bB_E\ \} 
    = \sup\{\|p(X)\| : X^T \in  \fC\bB_E\ \} \\
    &= \sup\{\|p(X)\| : X \in  \fC\bB_F\ \} = \|p\|_{A(\fC\bB_F)} . 
\end{align*}
Thus, the identity map extends to an isometric map $A(\fC\bB_E) \to A(\fC\bB_F)$. 
We can naturally identify these algebras of functions on $\widetilde{\fC\bB}_E = \widetilde{\fC\bB}_F$, since a commuting tuple $X$ is similar to a strict row/column contraction if and only if its joint spectrum is in the open ball, and therefore these similarity envelopes are indeed equal. 

However, the identity map $A(\fC\bB_E) \to A(\fC\bB_F)$ is not completely bounded. 
We shall use the fact that $A(\fC\bB_E)$ can be identified with $\cA_d$ --- the norm closure of polynomials in the algebra $\cM_d = \mlt(H^2_d)$ of multipliers on the Drury--Arveson space (in a sense this observation goes back, see e.g.  \cite{jury2012universal}; see also \cite[Section 11]{SSS1} or the section ``Connections to Noncommutative Function Theory" in \cite{hartz2025operator}). 
We shall also use the well-known fact that the multiplication operators $M = (M_1, 
\ldots, M_d) = (M_{z_1}, \ldots, M_{z_d})$ form a universal commuting row contraction, and hence $(M_1^*, \ldots, M_d^*)$ is a universal commuting column contraction. 
Therefore, to show that the identity map is not completely bounded, it suffices to find a sequence of matrix valued polynomials $P_n $ such that 
\[
\|P_n(M_1, \ldots, M_d)\| \leq 1 \,\, \textrm{ and } \,\, \sup_n\|P_n(M_1^*, \ldots, M_d^*)\| = \infty .
\]
To this end, consider the matrix valued polynomials
\[
P_n(z) = \operatorname{row}(z^w)_{|w| = n}
\]
that is, the row formed by all monomials $z^w = z_{w_1} \cdots z_{w_n}$ with words $w \in \{1, \ldots, d\}^n$ of length $n$ (note that we are using NC power notation and that most commutative monomials $z^\alpha = z_1^{\alpha_1} \cdots z_d^{\alpha_d}$ with $\alpha \in \bN^d$ will appear more than once since the variables are commutative). 
Then since $M$ is a row contraction, $P_n(M)$ is itself a row contraction of length $d^n$, in other words $\|P_n(M)\| \leq 1$. 

On the other hand, $\|P_n(M_1^*, \ldots, M_d^*)\|^2$ is the square of the row norm of $P_n(M_1^*, \ldots, M_d^*)$, which is equal to 
\begin{align*}
\left\|\sum_{|w|=n} (M^w)^* M^w \right\|
& \geq 
\left\langle \sum_{|w|=n} (M^w)^* M^w 1, 1 \right\rangle \\
&= \sum_{|w|=n} \|z^w\|^2 \\
&= \frac{(n+d-1)!}{n!(d-1)!}, 
\end{align*}
where we have used the well-known value of the norm sum (see \cite[Corollary 4.9]{Arv-subalg3}). 
It follows that $\|P_n(M_1^*, \ldots, M_d^*)\| \xrightarrow{n \to \infty} \infty$. This proves that the identity map is not completely bounded. 

In fact, there is no completely bounded isomorphism whatsoever from $A(\fC\bB_E)$ onto $A(\fC\bB_F)$; see Example \ref{ex:row_col_not_cb}. 
\end{example}

The following is a finite-dimensional example for which we do not spell out the details. 
\begin{example}
Consider the varieties
\[
\fV = \{(X,Y) \in \bB_{(\ell_2^2)_{\operatorname{row}}} : X^2 = XY = YX = Y^2 = 0\} 
\]
and
\[
\fW = \{(X,Y) \in \bB_{(\ell_2^2)_{\operatorname{col}}} : X^2 = XY = YX = Y^2 = 0\} .
\]
Then $A(\fV)$ and $A(\fW)$ are isometric finite-dimensional operator algebras which are not completely isometric. 
\end{example}

\subsection{NC maps need not extend off subvarieties}\label{subsec:maps_dont_extend}

In \cite[Theorem 6.5]{SampatShalit25}, a nonlinear injectivity phenomenon was put forward: if $G_0 \colon \fV \to \bB_F$ is an NC map from an NC subvariety $\fV \subseteq \bB_E$ into an NC operator ball $\bB_F$, and if $F$ is an injective operator space, then $G_0$ extends to an NC map $G \colon \bB_E \to \bB_F$. 
This result was used to sharpen the classification results by showing that certain NC function algebras are classified not by maps between varieties that carry one onto the other, but rather by maps between the balls that carry one variety onto the other. 

It is natural to seek a similar extension theorem in the spectral setting, that is, when balls and subvarieties are replaced by their similarity envelopes. 
The following example shows that such an extension does not hold. 

\begin{example}\label{ex:not_extend}
    Let $E = F = M_2$ be the operator space of $2 \times 2$ matrices, which is injective. 
    Let $\fV = \fC\bB_E$ be the block $2 \times 2$ matrices with commuting block entries. 
    We will show that there exists an NC map $G_0 \colon \widetilde{\fV} \to \widetilde{\bB}_F$  that does not extend to an NC map $G \colon \widetilde{\bB}_E \to \widetilde{\bB}_F$, in contrast with the situation with maps into balls described above.   
    Let $\tau \colon M_2 \to M_2$ be the transpose map $\tau(B) = B^T$, and consider the NC map $G_0 = \tau|_{ \widetilde{\fV}}$, that is, $G_0$ acts on $M_n(M_2)$ as $I_n \otimes \tau$, or on $M_2(M_n)$ as 
    \[
    G_0\left(\begin{pmatrix} X_{11} & X_{12} \\ X_{21} & X_{22} \end{pmatrix} \right) = \begin{pmatrix} X_{11} & X_{21} \\ X_{12} & X_{22} \end{pmatrix} .
    \]
    The fact that $G_0$ maps $\fV$ into $\widetilde{\bB}_F$ follows from \cite[Theorem 4.4]{SchShaSha}, which says that for commuting matrices $X_{11}, X_{12}, X_{21}, X_{22}$
    \[ 
    \rho_F \left(
    \begin{pmatrix}
    X_{11} & X_{12} \\
    X_{21} & X_{22}
    \end{pmatrix}
    \right) = \max \left\{\left\| \begin{pmatrix} \lambda_{11} & \lambda_{12} \\ \lambda_{21} & \lambda_{22}  \end{pmatrix} \right\| : \lambda \in \sigma(X) \right\},
    \]
    where $\sigma(X) = \sigma((X_{11} , X_{12}, X_{21}, X_{22}))$ is the joint spectrum of the commuting tuple $X = (X_{11} , X_{12}, X_{21}, X_{22})$.

    Let $G \colon \widetilde{\bB}_E \to \bM(M_2) \cong \bM^4$ be an NC holomorphic map that extends $G_0$. 
    Using the definition of the NC difference-differential operator $\Delta G(X,Y)$ to be the linear map that makes the following top-right block correct
    \[
    G \left(\begin{pmatrix} X & Z \\ 0 & Y \end{pmatrix} \right) = \begin{pmatrix} G(X) & \Delta G(X,Y) [Z] \\ 0 & G(Y) \end{pmatrix}, 
    \]
    for $X \in M_m(M_2), Y \in M_n(M_2)$ and $Z \in M_{m,n}(M_2)$, we find that $\Delta G(0,0)[Z] = Z^T = \tau(Z)$, where $0 = 0_2 \in M_2$ and $Z \in M_2$. 
    By the properties of NC maps it follows that if $0$ denotes the $n$-fold ampliation of the scalar $0$, then $\Delta G(0,0) = I_n \otimes \tau$ is the ampliation of the transpose, for all $n$. 
    Loosely, we write this as $\Delta G(0,0) = \tau$. 

    But now if we take 
    \[
    X = \begin{pmatrix} E_{11} & E_{21} \\ E_{12} & E_{22} \end{pmatrix}, 
    \]
    then $X$ is unitary, therefore $\frac{3}{4} X \in \bB_E$. 
    On the other hand, 
    \[
    \tau\left(\frac{3}{4}X\right) = \tau\left(\frac{3}{4} \begin{pmatrix} E_{11} & E_{21} \\ E_{12} & E_{22} \end{pmatrix} \right) = \frac{3}{4} \begin{pmatrix} E_{11} & E_{12} \\ E_{21} & E_{22} \end{pmatrix}.
    \]
    The maximal eigenvalue of the right hand side is $\frac{3}{2}$.  
    Therefore, $\tau(\frac{3}{4}X)$ is not similar to a contraction, so in particular there does not exist an invertible $S \in M_2$ such that $[I_2 \otimes S^{-1}]\tau(\frac{3}{4}X) [I_2 \otimes S] \in \bB_F$, and so $\rho_E(\frac{3}{4}X) < 1$ while $\rho_F(\tau(\frac{3}{4}X)) > 1$. 
    
    On the other hand, it follows from Lemma \ref{lem:spectral_less} that if $G$ maps $\widetilde{\bB}_E$ to $\widetilde{\bB}_F$, then $\rho_F(\Delta G(0,0)[Z]) \leq \rho_E(Z)$, thus 
    \[
    \rho_F(\tau(Z)) \leq \rho_E(Z)
    \]
    for every $Z$. 
    We conclude that any NC holomorphic extension $G  \colon \widetilde{\bB}_E \to \bM(M_2)$ of $G_0$ cannot map into $\widetilde{\bB}_F$. 

    We note that by using the same ideas as in Example \ref{ex:bdd_not_cb}, it can be shown that the map $f \mapsto f \circ \tau$ gives rise to an  isometric isomorphism of $A(\fC \bB_E)$ onto itself, which is not completely bounded. 
\end{example}

\begin{remark}
Examples \ref{ex:bdd_not_cb} and \ref{ex:not_extend} show that for homogeneous NC subvarieties, there exist isometric isomorphisms which are not completely bounded. 
This is in contrast with the case of the row ball, where an isometric isomorphism must be completely isometric, and in fact must be induced by a restriction of an automorphism of the ambient operator ball \cite{SSS1,SSS2}. 
These examples also present an interesting contrast to the completely isometric classification of NC function algebras on subvarieties of arbitrary NC operator balls \cite{SampatShalit25,sampat2025weak} --- under reasonable assumptions, a completely isometric isomorphism $\alpha \colon A(\fV) \to A(\fW)$ for homogeneous varieties $\fV \subseteq \bB_E$ and $\fW \subseteq \bB_F$ arises from an isomorphism between the two NC balls, implying that $E$ and $F$ are completely isometric (see Theorems 6.12 and 6.21 in \cite{SampatShalit25}). 
\end{remark}

\section{Rigidity of NC ball algebras}

In this section we consider the isomorphism problem for algebras $A(\fV)$ in the extreme case where the variety $\fV$ consists of the entire NC unit ball $\bB_E$ for some operator space $E$. 
We will see that NC ball algebras $A(\bB_E)$ and $A(\bB_F)$ can be completely boundedly isomorphic only when $E$ and $F$ are completely isometric, and hence the algebras are already completely isometrically isomorphic. 

It will be useful to recall that over every operator space $\cE$, one may form the universal unital operator algebra $OA_u(\cE)$ over $\cE$ (see \cite[pp. 112--113]{Pisier-book}).
By \cite[Section 7.1]{SampatShalit25}, the algebra $A(\bB_E)$ can be naturally identified completely isometrically with $OA_u(\cE)$ where $\cE = E^*$ can be identified with $\spn\{Z_1, \ldots, Z_d\}$ for the coordinate functions $Z_1, \ldots, Z_d$ in $A(\bB_E)$, and $n$th degree monomials correspond to tensors in the Haagerup tensor power $\cE \otimes_h \cdots  \otimes_h \cE$ of $\cE$ with itself $n$ times.  

\begin{theorem}\label{thm:ball_alg_iso}
Let $E$ and $F$ be finite-dimensional operator spaces and let $\bB_E$ and $\bB_F$ be the corresponding NC operator unit balls. 
Suppose that there exists a completely bounded isomorphism $\alpha \colon A(\bB_E) \to A(\bB_F)$. 
Then $E$ and $F$ are completely isometric.
\end{theorem}
\begin{proof}
First, observe that by Theorem \ref{thm:iso_vacuum_iso}, we may assume that $\alpha^* \Phi_0 = \Phi_0$. 
This means that $\alpha$ maps $A(\bB_E)_0 := \ker \Phi_0$ into $A(\bB_F)_0$. 
Note that $(\alpha^{-1})^* \Phi_0 = \Phi_0$ too, so that $\alpha$ maps $A(\bB_E)_0$ onto $A(\bB_F)_0$.

Let $P_n$ be the projection \eqref{eq:homog_proj} onto the $n$th homogeneous component. 
Define $\cE_n := \spn\{Z^w : |w|=n\} = P_n A(\bB_E)$ and likewise $\cF_n := P_n A(\bB_F)$. 
Thus, by the discussion before the theorem, $\cE_1 = E^*$ and $\cF_1 = F^*$. 
Define also 
\[
\alpha_n \colon \cE_n \to \cF_n \quad , \quad \alpha_n = P_n \alpha \big|_{\cE_n}.
\]
We shall prove that $\alpha_1 \colon \cE_1 \to \cF_1$ is a completely isometric isomorphism. 
Since $\cE_1^* = E^{**} = E$ and $\cF_1^* = F^{**} = F$, this will complete the proof. 

Clearly, $\alpha_n$ is completely bounded and $\|\alpha_n\|_{cb} \leq \|\alpha\|_{cb}$ for all $n$. 
Identifying $\cE_n$ with the Haagerup tensor product $\cE_1 \otimes_h \cdots \otimes_h \cE_1$, and using the fact that $\alpha$ is an isomorphism, we see that $\alpha_n = \alpha_1 \otimes_h \cdots \otimes_h \alpha_1$. 
By properties of the Haagerup tensor product we obtain 
\[
\|\alpha_1\|_{cb}^n = \|\alpha_n\|_{cb} \leq \|\alpha\|_{cb} , 
\]
for all $n \in \bN$. It follows that $\alpha_1$ is completely contractive. 

But because $\alpha$ is an algebraic isomorphism that maps $A(\bB_E)_0$ onto $A(\bB_F)_0$, we must have that $\alpha$ preserves the minimal degree of a monomial in the Taylor series in \ref{eq:homog_exp}. Hence $\alpha_1$ is injective. Furthermore, 
\[
P_1\alpha(e)=P_1\alpha P_1(e)+P_1\alpha(e-P_1(e))=P_1\alpha P_1(e)
\]
for all $e\in A(\bB_E)_0$. Thus for $f\in A(\bB_F)_0$, using that $\alpha^{-1}$ is vacuum preserving, we see that $\alpha^{-1}(f) \in A(\bB_E)_0$, so that 
\[
\alpha_1P_1\alpha^{-1}(f)=P_1\alpha P_1\alpha^{-1}(f)=P_1\alpha\alpha^{-1}(f)=P_1(f),
\]
showing that $\alpha_1$ is surjective. Moreover this equality shows that
\[
(\alpha^{-1})_1(f)=P_1\alpha^{-1}P_1(f)=P_1(\alpha_1)^{-1}\alpha_1P_1\alpha^{-1}P_1(f)=P_1(\alpha_1)^{-1}P_1(f)=(\alpha_1)^{-1}(f)
\]
for all $f\in\cF_1$. Thus $(\alpha_1)^{-1} = (\alpha^{-1})_1$.

By the same argument as above $(\alpha_1)^{-1}$ is completely contractive too. 
In conclusion, $\alpha_1\colon \cE_1 \to \cF_1$ is a completely isometric isomorphism, as required. 
\end{proof}

\begin{remark}
    Note that if the isomorphism $\alpha$ is only bounded, then the proof shows that the operator spaces are isometrically isomorphic. This follows also from Theorem \ref{thm:classification_linear}. 
\end{remark}
\begin{remark}
    It is interesting to compare Theorem \ref{thm:ball_alg_iso} to the recent \cite[Corollary 4.1]{tornes2026}, which says that if $A(\bB_E) = A(\bB_F)$ are equal with the identity a completely bounded map, then the identity map $E \to F$ is a complete isometry. 
\end{remark}

\begin{example}
    Let $E = \left(\ell^2_d\right)_{\operatorname{row}}$ and $F = \left(\ell^2_d\right)_{\operatorname{col}}$ be the row and column spaces over $\bC^d$ for some $d \geq 2$.  
    Then since $E$ and $F$ are not completely isometrically isomorphic, Theorem \ref{thm:ball_alg_iso} says that there can be no completely bounded isomorphism between the ball algebras $A(\bB_E)$ and $A(\bB_F)$, in particular the identity map sending NC polynomials from one algebra to the other is not completely bounded. 
    In \cite[Remark 8.5]{SchShaSha} it is shown that the identity map is not even bounded on NC polynomials (and therefore $A(\bB_{E}) \neq A(\bB_F)$ as sets), but this does not follow from the above remark because $E$ and $F$ are isometric as normed spaces. 
\end{example}

\begin{remark}
In the case that $E$ and $F$ are not finite dimensional, the proof of Theorem \ref{thm:ball_alg_iso} shows that if $A(\bB_E)$ and $A(\bB_F)$ are completely boundedly isomorphic via a vacuum preserving isomorphism, then $E^*$ and $F^*$ are completely isometrically isomorphic (note that Theorem \ref{thm:iso_vacuum_iso}, which allowed us to reduce to the vacuum preserving case, was proved under a standing finite dimensionality assumption). 
The non-uniqueness of the operator space predual implies that the conclusion is false in general, that is, there exist non-completely-isometric operator spaces $E_i$ such that $A(\bB_{E_1})$ is completely isometrically isomorphic to $A(\bB_{E_2})$ via the identity map (which is certainly vacuum preserving). 
This observation also holds in the bounded/isometric categories. 
\end{remark}

Since $A(\bB_E) = OA_u(\cE)$ for $\cE = E^*$, Theorem \ref{thm:ball_alg_iso} is equivalent to the following corollary, whose equivalent formulation may be of interest in itself. 
Although, by universality, it is clear that if $\cE$ and $\cF$ are completely isometric then $OA_u(\cE)$ and $OA_u(\cF)$ are completely isometrically isomorphic, the converse implication is not immediate, and as far as we know was first observed in \cite[Section 7.1]{SampatShalit25}. 
The corollary implies something quite stronger: already the complete isomorphism class of $OA_u(\cE)$ determines $\cE$ up to complete isometry. 

\begin{corollary}\label{cor:univ_alg_iso}
Let $\cE$ and $\cF$ be finite-dimensional operator spaces. 
If the universal operator algebras $OA_u(\cE)$ and $OA_u(\cF)$ are completely boundedly isomorphic, then $\cE$ and $\cF$ are completely isometric.
\end{corollary}

\section{On the classification of commutative ball algebras}\label{sec:comm_balls}

Following the strong rigidity results of the previous section, it is natural to ask what about the commutative slices of the NC balls? Let $\fV = \fC\bB_E$ and $\fW = \fC\bB_F$. Of course our methods show that if $A(\fV)$ and $A(\fW)$ are boundedly isomorphic then $E$ and $F$ are isometric, hence $E$ and $F$ are completely boundedly isomorphic by finite dimensionality. What about the other implications? 
If $A(\fV)$ and $A(\fW)$ are completely boundedly isomorphic, does it follow that $E$ and $F$ are completely isometric? 
The following example shows that the answer is no. 

\begin{example}\label{ex:cb_not_ciso}

Let $E=\left(\ell^2_2\right)_{\operatorname{row}}$, and let $F$ be the operator space obtained from the embedding
\[
\mathbb C^2\longrightarrow \left(\ell^2_2\right)_{\operatorname{row}}\oplus_\infty \left(\ell^2_2\right)_{\operatorname{col}},
\qquad
x\longmapsto \left(x,\frac34x\right).
\]
For $X=(X_1,X_2)\in M_n(\mathbb C^2)$,
\[
\|X\|_{M_n(F)}
=
\max\left\{
\left\|X_1X_1^*+X_2X_2^*\right\|^{1/2},
\frac34\left\|X_1^*X_1+X_2^*X_2\right\|^{1/2}
\right\}.
\]
Put
\[
\fV=\fC\bB_E 
\quad \textrm{ and } 
\quad
\fW=\fC\bB_F.
\]
Then $E$ and $F$ are isometric (since at the scalar level, both $E$ and $F$ have the Euclidean norm) but not completely isometric. 
Indeed, if
\[
X_1=E_{11},\qquad X_2=E_{21},
\]
then
\[
\|(X_1,X_2)\|_{M_2(E)}=1,
\qquad
\|(X_1,X_2)\|_{M_2(F)}
=\frac{3\sqrt2}{4}>1.
\]
Since every linear isometry of $\ell_2^2$ is unitary, and coordinate unitaries preserve both the row and column norms, $E$ and $F$ cannot be completely isometric.

We shall see that $A(\fV)$ is completely boundedly isomorphic to $A(\fW)$.
Thus, the function algebras on commutative slices of NC operator balls do not enjoy the same kind of rigidity that the function algebras on the full NC balls enjoy.

We now compare the two function algebras. Since $\fW\subseteq\fV$, every matrix-valued polynomial $P$ satisfies
\[
\|P\|_{A(\fW)}
\leq
\|P\|_{A(\fV)}.
\]
Let $H^2_d$ be the Drury--Arveson space and let $M=(M_1,M_2)$ be the shift (as in Example \ref{ex:bdd_not_cb}). 
Since $M$ is the universal commuting row contraction,
\[
\|P\|_{A(\fV)}=\|P(M)\|.
\]
Write $H^2_d$ as the direct sum $\oplus_{n\geq 0} H_n$ of spaces $H_n$ of homogeneous polynomials of degree $n$. 
Set
\[
q=\frac{2\sqrt2}{3},
\qquad
D|_{H_0}=I,
\qquad
D|_{H_n}=qI\quad(n\geq1),
\]
and define $T_i=DM_iD^{-1}$. 
Direct computations show that 
\[
T_1 T_1^* + T_2 T_2^*
=
0_{H_0}\oplus q^2I_{H_1}
\oplus I_{\bigoplus_{n\geq2}H_n}
\leq I
\]
and
\[
T_1^* T_1 + T_2^* T_2
=
\frac{16}{9}I_{H_0}
\oplus
\bigoplus_{n\geq1}\frac{n+2}{n+1}I_{H_n}
\leq \frac{16}{9}I.
\]
Thus, every compression of $T$ to a finite-dimensional subspace has $\|\cdot\|_F$ norm at most $1$. By compressing $T$ to
$\bigoplus_{n=0}^N H_n$, and then applying a radial scaling, one obtains sufficiently many commuting tuples in $\fW$ to approximate $T$; consequently,
\[
\|P(T)\|\leq \|P\|_{A(\fW)}.
\]
Since $M=D^{-1}TD$, it follows, at every matrix level, that
\[
\begin{aligned}
\|P\|_{A(\fV)}
&=\|P(M)\| \\
&\leq \|D^{-1}\|\,\|D\|\,\|P(T)\| \\
&\leq \frac{3\sqrt2}{4}\,
       \|P\|_{A(\fW)}.
\end{aligned}
\]
Therefore
\[
\|P\|_{A(\fW)}
\leq
\|P\|_{A(\fV)}
\leq
\frac{3\sqrt2}{4}\,
\|P\|_{A(\fW)}
\]
completely, and the identity on polynomials extends to a completely bounded isomorphism
$A(\fV)\to A(\fW)$, as claimed.
\end{example}

The following example complements the example above --- it exhibits two naturally occurring operator spaces for which the corresponding commutative ball algebras are different. 
Before finding the above example, this might have served as evidence supporting the conjecture that commutative ball algebras are rigid. 
In light of the previous example, the main role of the following one is to illustrate how our methods can be used to separate algebras.

\begin{example}\label{ex:row_col_not_cb}
Let $E$ be the row operator space, let $F$ be the column operator space in $d \geq 2$ dimensions, and let $\fC \bB_E$ and $\fC \bB_F$ be the respective commutative balls, as in Example \ref{ex:bdd_not_cb}. 
In Example \ref{ex:bdd_not_cb} we showed that the identity map on polynomials does not extend to a completely bounded isomorphism $A(\fC\bB_E)$ onto $A(\fC\bB_F)$. 
We shall now show that there is no completely bounded isomorphism whatsoever between these two algebras. 

Indeed, suppose that there is a completely bounded isomorphism of $A(\fC\bB_E)$ onto $A(\fC\bB_F)$. 
By Theorem \ref{thm:classification_linear}, there exists a linear bijection $L \colon F \to E$ that promotes to a $\delta_{cb}$-bi-Lipschitz map between the similarity envelopes $\widetilde{\fC\bB}_F$ and $\widetilde{\fC\bB}_E$. 
Then $L$ itself induces, by Proposition \ref{prop:equiv_G_norm}, a completely bounded isomorphism $\alpha$ between the algebras. 
Since the first level of both envelopes is $\bB_d$, $L$ must be a unitary map. But the unitary map $L^{-1}$ induces an automorphism $\beta$ of $A(\fC\bB_E)$. 
Composing the automorphism $\beta$ with the isomorphism $\alpha$, we obtain a completely bounded isomorphism $\alpha \circ \beta$ of $A(\fC\bB_E)$ onto $A(\fC\bB_F)$ which is induced by composition with the identity. 
This shows that the identity mapping is a completely bounded isomorphism of $A(\fC\bB_E)$ onto $A(\fC\bB_F)$, in contradiction to what we showed in Example \ref{ex:bdd_not_cb}. 
\end{example}

\section{Classification of algebras of bounded NC functions}

In this section, we will briefly describe the classification problem for the algebras of the type $H^\infty(\fV)$, where $\fV$ is an arbitrary (that is, not necessarily homogeneous) subvariety\footnote{Recall that a variety in $\bB_E$ is the joint zero set of a family of $H^\infty(\bB_E)$ functions.} of $\bB_E$ for a finite-dimensional operator space structure $E$. 
We shall classify these algebras up to bounded and completely bounded weak-$*$ continuous isomorphism. 
The role of the weak-$*$ continuity assumption is that, thanks to Theorem \ref{thm:weakstar_rho_lessthan_1}, it guarantees that an isomorphism $\alpha$ induces an adjoint $\alpha^*$ which maps a point evaluation $\Phi_W$ to a point evaluation $\Phi_{G(W)} :=\alpha^* \Phi_W$, thus induces a map $G$ between the similarity envelopes of the varieties. 
One might conjecture that isomorphisms map point evaluations to point evaluations automatically, but Example \ref{ex:non_singleton_fiber} presents an obstruction to proving this by using function theoretic means to keep track of the points over which a representation is fibered. 

The following lemma shows that the natural assumption of weak-$*$ continuity is equivalent to the operationally convenient assumption of preserving point evaluations.

\begin{lemma}\label{lem:weak_point}
Let $\alpha \colon H^\infty(\fV) \to H^\infty(\fW)$ be a unital homomorphism. 
Then $\alpha$ is weak-$*$ continuous if and only if $\Phi_W \circ \alpha$ is a point evaluation for every point evaluation $\Phi_W,\, W \in \widetilde{\fW} $. 
\end{lemma}
\begin{proof}
By Theorem \ref{thm:weakstar_rho_lessthan_1}, a finite-dimensional representation $\Phi \in \operatorname{Rep}_{\mathrm{fin}}(H^\infty(\fW))$ is weak-$*$ continuous if and only if $\Phi$ is a point evaluation $\Phi_W$ for some $W \in \widetilde{\fW}$. 
If $\alpha$ is weak-$*$ continuous, then $\alpha^*(\Phi_W) = \Phi_W\circ \alpha$ is weak-$*$ continuous for all $W \in \widetilde{\fW}$ and is therefore a point evaluation. 
It follows that $\alpha^*$ maps point evaluations to point evaluations. 

Conversely, suppose that $\Phi_W \circ \alpha$ is a point evaluation, say $\Phi_{G(W)}$ for every point evaluation $\Phi_W$ (here $W \in \widetilde{\fW}$ and $G(W) \in \widetilde{\fV}$). 
We first claim that $\alpha$ must be bounded. 
To see this, note that for every $f \in H^\infty(\fV)$ we have that 
\[
\sup_{W \in \fW} \|\Phi_W\circ\alpha(f)\| = \sup_{W \in \fW} \|\alpha(f)(W)\| < \infty, 
\]
because $\alpha(f) \in H^\infty(\fW)$, by assumption. 
By the uniform boundedness principle, 
\[
C : = \sup_{W \in \fW} \|\Phi_W\circ\alpha\| < \infty, 
\]
leading to $\|\alpha(f)\| = \sup_{W\in \fW} \|\alpha(f)(W)\|  \leq C\|f\|$ for all $f \in H^\infty(\fV)$. 
Thus $\alpha$ is bounded. 

Now if $f_n$ is a bounded net in $H^\infty(\fV)$ converging weak-$*$ to $f$, then for every $W \in \fW$, 
\begin{align*}
\alpha(f_n)(W) 
&= \Phi_W \circ \alpha(f_n) \\
&= \Phi_{G(W)} (f_n) \to \Phi_{G(W)}(f) = \alpha(f)(W), 
\end{align*}
and this shows that $\alpha$ takes a bounded weak-$*$ convergent net to a bounded weak-$*$ convergent net. The map $\alpha$ is therefore weak-$*$ continuous on the closed unit ball, and so by the Krein-Smulian theorem $\alpha$ is weak-$*$ continuous, as required. 
\end{proof}

\subsection{Maps induce homomorphisms}

The following proposition is an analogue of Proposition \ref{prop:equiv_G_norm}. Note the difference that we do not need to assume that $G\big|_\fW$ is bounded, since this follows from each of the equivalent conditions. 
We also require no homogeneity assumption.

\begin{proposition}\label{prop:equiv_G_infty}
Let $\fV \subseteq \bB_E$ and $\fW \subseteq \bB_F$ be NC varieties, and $G \colon \widetilde{\fW} \to \widetilde{\fV}$ an NC holomorphic function. Let $\kappa \in \{b, cb\}$. Then the following statements are equivalent:
\begin{enumerate}
\item $G$ is Lipschitz with respect to $\delta_{\kappa}$,
\item $\sup_{W\in \fW}\|\Phi_{G(W)}\|_{\kappa} < \infty$, ~ and
\item the mapping $f \mapsto \widetilde{f} \circ G$ for $f \in \cH^\infty(\fV)$ is a well-defined homomorphism $\alpha:H^\infty(\fV) \to H^\infty(\fW)$, which is completely bounded in the case $\kappa = cb$.
\end{enumerate}
Furthermore, in this case
\begin{itemize}
\item[(a)] $\|\alpha \|_{\kappa} = \sup_{W\in \fW}\|\Phi_{G(W)}\|_{\kappa}$,
\item[(b)] the value in (a) is a Lipschitz constant for $G$, ~ and
\item[(c)] $\alpha$ is weak-$*$ continuous.
\end{itemize}
\end{proposition}
\begin{proof}
Most of the proof is similar to that of Proposition \ref{prop:equiv_G_norm}, so we omit it. 
Note that Lemma \ref{lem:weak_point} gives item $(c)$. 
We shall give some details for the implication $(1)\Longrightarrow(2)$ since the proof of Proposition \ref{prop:equiv_G_norm} made use of homogeneity in that step. 

Let $L$ be a $\delta_\kappa$-Lipschitz constant for $G$, 
and fix some 
$W_0\in\fW(m) \subset \fW,$ for some $m\geq1$. 
Let $W\in\fW(n)$. Then letting $W^{(m)} = W \oplus \cdots \oplus W$ ($m$ times), we have that $\|\Phi_{G(W^{(m)})}\|_\kappa = \|\Phi_{G(W)}\|_\kappa$. 
Moreover, 
\begin{align*}
\|\Phi_{G(W^{(m)})}\|_\kappa
&\leq
\|\Phi_{G(W^{(m)})}-\Phi_{G(W_0^{(n)})}\|_\kappa
+
\|\Phi_{G(W_0^{(n)})}\|_\kappa \\
&=
\delta_\kappa\bigl(G(W^{(m)}),G(W_0^{(n)})\bigr)
+
\|\Phi_{G(W_0)^{(n)}}\|_\kappa \\
& \leq
L\,\delta_\kappa(W^{(m)},W_0^{(n)})
+
\|\Phi_{G(W_0)}\|_\kappa.
\end{align*}
Since the evaluation representations
$\Phi_W$ and $\Phi_{W_0}$ are completely contractive, $\delta_\kappa(W^{(m)},W_0^{(n)}) \leq 2$, we conclude that 
\[
\|\Phi_{G(W^{(m)})}\|_\kappa \leq 2L + \|\Phi_{G(W_0)}\|_\kappa  < \infty, 
\]
and this shows $(2)$. 
\end{proof}

\subsection{Homomorphisms induce maps}

\begin{proposition}\label{prop:weakstar_homo_ind_NCmap}
    Let $\fV \subseteq \bB_E$ and $\fW \subseteq \bB_F$ be varieties, and let $\alpha \colon H^\infty(\fV) \to H^\infty(\fW)$ be a weak-$*$ homeomorphic isomorphism. Then there exists an NC biholomorphic map $G \colon \widetilde{\fW} \to \widetilde{\fV}$ that is bi-Lipschitz with respect to $\delta_b$ ($\delta_{cb}$ if $\alpha$ is a completely bounded isomorphism) and $\alpha$ is given by 
    \[
    \alpha(f) = \widetilde{f}\circ G. 
    \]
\end{proposition}
\begin{proof}
    By Theorem \ref{thm:weakstar_rho_lessthan_1}, we can identify the finite-dimensional weak-$*$ continuous representations of $H^\infty(\fV)$ with the point evaluations $\Phi_X$ at points $X \in \widetilde{\fV}$. 
    Arguing as in the proof of Proposition \ref{prop:homo_ind_NCmap} gives the result.
\end{proof}

\subsection{Classification}

\begin{corollary}\label{cor:Hinf_class}
    Let $\fV \subseteq \bB_E$ and $\fW \subseteq \bB_F$ be varieties. 
    There exists a weak-$*$ homeomorphic isomorphism $H^\infty(\fV) \to H^\infty(\fW)$ if and only if $\widetilde{\fV}$ and $\widetilde{\fW}$ are NC biholomorphic via a $\delta_b$-bi-Lipschitz map $G$. There exists a weak-$*$ homeomorphic and completely bounded isomorphism $H^\infty(\fV) \to H^\infty(\fW)$ if and only if $\widetilde{\fV}$ and $\widetilde{\fW}$ are NC biholomorphic via a $\delta_{cb}$-bi-Lipschitz map.  
\end{corollary}

We obtain a bit more in the homogeneous case. 

\begin{theorem}\label{thm:classification_linear_weak_star}
Let $\fV \subseteq \bB_E$ and $\fW \subseteq \bB_F$ be homogeneous varieties in finite-dimensional operator balls and let $\kappa \in \{b,cb\}$. 
Then there exists a weak-$*$ homeomorphic $\kappa$-bounded isomorphism $H^\infty(\fV) \to H^\infty(\fW)$ if and only if there exist linear maps $L \colon F \to E$ and $M \colon E \to F$ such that $L$ promotes to a bijective $\delta_\kappa$-bi-Lipschitz map from $\widetilde{\fW}$ onto $\widetilde{\fV}$ and $M$ promotes to its inverse $\widetilde{\fV} \to \widetilde{\fW}$. 
\end{theorem}
\begin{proof}
    The backward implication follows from Corollary \ref{cor:Hinf_class}. For the forward implication, assume that there is a weak-$*$ homeomorphic $\kappa$-bounded isomorphism $H^\infty(\fV) \to H^\infty(\fW)$. 
    By Corollary \ref{cor:Hinf_class}, there is $\delta_\kappa$-bi-Lipschitz map from $\widetilde{\fW}$ onto $\widetilde{\fV}$. 
    By Theorem \ref{thm:linearisation_maps}, there is a $\kappa$-bounded isomorphism $A(\fV) \to A(\fW)$. 
    By Theorem \ref{thm:classification_linear}, we obtain the linear maps $L \colon F \to E$ and $M \colon E \to F$ that restrict to $\delta_\kappa$-bi-Lipschitz bijections between $\widetilde{\fW}$ and $\widetilde{\fV}$, as required. 
\end{proof}

By the previous theorem together with Theorem \ref{thm:classification_linear}, we see that for homogeneous varieties $H^\infty(\fV)$ is weak-$*$ homeomorphically and (completely) boundedly isomorphic to $H^\infty(\fW)$, if and only if $A(\fV)$ is (completely) boundedly isomorphic to $A(\fW)$. 
Therefore, one can see that the examples of Section \ref{sec:comm_balls} persist in the $H^\infty$ setting. 
Likewise, an immediate application of Theorem \ref{thm:ball_alg_iso} gives the following corollary.

\begin{corollary}\label{cor:ball_alg_isoHinf}
Let $E$ and $F$ be finite-dimensional operator spaces and let $\bB_E$ and $\bB_F$ be the corresponding NC operator unit balls. 
Suppose that there exists a weak-$*$ homeomorphic completely bounded   $\alpha \colon H^\infty(\bB_E) \to H^\infty(\bB_F)$. 
Then $E$ and $F$ are completely isometric.
\end{corollary}

\bibliographystyle{abbrv}
\bibliography{bibliography}

\end{document}